\documentclass[11pt, reqno]{amsart}
\usepackage{amsfonts,latexsym}
\usepackage{amsmath}
\usepackage{amscd}
\usepackage{float,amsmath,amssymb,mathrsfs,bm,multirow,graphics}
\usepackage[dvips]{graphicx}
\usepackage[percent]{overpic}
\usepackage[numbers,sort&compress]{natbib}
\usepackage{enumerate}
\usepackage{amsaddr}

\newcommand{\nequation}{\setcounter{equation}{0}}

\newcommand{\R}{{\Bbb R}}

\newcommand{\C}{{\Bbb C}}

\newcommand{\Z}{{\Bbb Z}}
\newcommand{\proofbegin}{\noindent{\it Proof.\,\,}}
\newcommand{\proofend}{\hfill$\Box$\bigskip}

\DeclareMathOperator{\tr}{tr}
\DeclareMathOperator{\sgn}{sgn}
\DeclareMathOperator{\im}{Im}
\DeclareMathOperator{\re}{Re}

\newcommand{\res}{\text{\upshape Res\,}}

\def\XXint#1#2#3{{\setbox0=\hbox{$#1{#2#3}{\int}$}
\vcenter{\hbox{$#2#3$}}\kern-.5\wd0}}

\allowdisplaybreaks

\newtheorem{theorem}{Theorem}[section]
\newtheorem{proposition}[theorem]{Proposition}

\newtheorem{definition}[theorem]{Definition}

\newtheorem{remark}[theorem]{Remark}
\newtheorem{example}[theorem]{Example}

\input epsf
\title[The nonlinear Schr\"odinger equation with $t$-periodic data II]{\sc The nonlinear Schr\"odinger equation with \\
$t$-periodic data: II. Perturbative results}

\author{J. Lenells}
\address{Department of Mathematics, KTH Royal Institute of Technology, \\ 100 44 Stockholm, Sweden.}
\email{jlenells@kth.se}

\author{A. S. Fokas}
\address{Department of Applied Mathematics and Theoretical Physics, University of Cambridge, Cambridge CB3 0WA, United Kingdom, and Research Center of Mathematics, Academy of Athens, 11527, Greece.}
\email{T.Fokas@damtp.cam.ac.uk} 
    
\begin{document}

\begin{abstract} 
\noindent
We consider the nonlinear Schr\"odinger equation on the half-line with a given Dirichlet boundary datum which for large $t$ tends to a periodic function. We assume that this function is sufficiently small, namely that it can be expressed in the form $\alpha g_0^b(t)$, where $\alpha$ is a small constant. Assuming that the Neumann boundary value tends for large $t$ to the periodic function $g_1^b(t)$, we show that $g_1^b(t)$ can be expressed in terms of a perturbation series in $\alpha$ which can be constructed explicitly to any desired order. As an illustration, we compute $g_1^b(t)$ to order $\alpha^8$ for the particular case that $g_0^b(t)$ is the sum of two exponentials. We also show that there exist particular functions $g_0^b(t)$ for which the above series can be summed up, and therefore for these functions $g_1^b(t)$ can be obtained in closed form. The simplest such function is $\exp(i\omega t)$, where $\omega$ is a real constant.
\end{abstract}

\maketitle

\noindent
{\small{\sc AMS Subject Classification (2010)}: 35Q55, 37K15.}

\noindent
{\small{\sc Keywords}: Initial-boundary value problem, time-periodic data, long-time asymptotics.}


\section{Introduction}\nequation
The problem of determining the Dirichlet to Neumann map for elliptic PDEs is well known. Actually, an analogous problem also exists for evolution PDEs. For eample, for the nonlinear Schr\"odinger (NLS) equation formulated on the half-line with given initial and boundary data, this problem amounts to expressing the spatial derivative at the boundary in terms of the given data. This problem is analyzed in \cite{trilogy1} and \cite{trilogy2} using two different formulations, both of which are based on the analysis of the so-called global relation: The formulation in \cite{trilogy1} is based on the eigenfunctions involved in the definition of the spectral functions $\{A(k), B(k)\}$ (see also \cite{DMS2001, F2005}), whereas the formulation in \cite{trilogy2} is based on an extension of the Gelfand-Levitan-Marchenko approach first introduced in \cite{BFS2003}.
It was shown in \cite{trilogy2} for the NLS, and in \cite{HF2013, H2014} for the mKdV and sine-Gordon equations, that if $u(0,t) = \alpha \sin{t}$, $\alpha \in \R$, then the function $u_x(0,t)$ for the NLS, and the functions $\{u_x(0,t), u_{xx}(0,t)\}$ for the mKdV and the sine-Gordon, can be computed {\it explicitly} at least up to and including terms of $O(\alpha^3)$, and furthermore the above functions become {\it periodic} as $t \to \infty$. 
Unfortunately, the perturbative approach of \cite{trilogy2} is quite cumbersome and it is practically impossible to go beyond terms of $O(\alpha^3)$.

Here we consider the NLS equation
\begin{align}\label{nls}
  iu_t + u_{xx} - 2\lambda |u|^2 u = 0, \qquad x>0, \quad t>0, \quad \lambda = \pm 1,
\end{align}
on the half-line and denote by $u_0(x)$ and $g_0(t)$  the given initial datum and the given Dirichlet boundary datum; we also denote by $g_1(t)$ the unknown Neumann boundary value:
\begin{align}\nonumber
& u_0(x) = u(x,0), && 0 < x < \infty; 
	\\ \label{star}
& g_0(t) = u(0,t), \qquad g_1(t) = u_x(0,t), && 0 < t < \infty.
\end{align}
We assume that $u_0(x) \in \mathcal{S}([0,\infty))$, where $\mathcal{S}([0,\infty))$
denotes the Schwartz class
\begin{align}\label{schwartzdef}
\mathcal{S}([0,\infty)) = \{u \in C^\infty([0,\infty)) \, | \, x^n u^{(m)}(x) \in L^\infty([0,\infty)) \text{ for all } n, m \geq 0\}.
\end{align}
Furthermore, we assume that $g_0(t)$ is asymptotically periodic as $t \to \infty$, namely,
\begin{align}\label{twostar}
  g_0(t) - g_0^b(t) = O(t^{-7/2}), \qquad t \to \infty, 
\end{align}
where $g_0^b(t)$ is a given periodic function of period $\tau > 0$.

We consider perturbative solutions of the NLS equation (\ref{nls}) with the initial and boundary conditions specified by equations (\ref{star}) and (\ref{twostar}). Under the basic assumption that there exists a solution for which 
\begin{align}\label{g1minusg1B}
g_1(t) - g_1^b(t) = O(t^{-7/2}), \qquad t\to \infty,
\end{align}
where $g_1^b(t)$ is a periodic function with period $\tau$, we show that $g_1^b(t)$ is uniquely determined from $g_0^b(t)$ to {\it all} orders in a perturbative expansion. Moreover, we characterize the coefficients of the Fourier series of $g_1^b(t)$ in terms of an infinite system of {\it algebraic} equations. It is straightforward to solve this algebraic system perturbatively and we show that the perturbative solution can be continued to all orders. Thus, under the basic assumption (\ref{g1minusg1B}) of the asymptotic periodicity of the Neumann boundary value, {\it this result provides an explicit construction of the Dirichlet to Neumann map for asymptotically $t$-periodic data in the limit of large $t$ and small data.}

As illustrations, we consider the example of a single exponential,
\begin{align}\label{1.1}
  g_0^b(t) = \alpha e^{i\omega t}, \qquad \alpha > 0, \quad \omega \in \R, \quad t > 0,
\end{align}
as well as the example of the sum of two exponentials
\begin{align}\label{1.5}
  g_0^b(t) = \alpha e^{i\omega t} + \beta e^{-i\omega t}, \qquad \alpha, \beta \in \C, \quad \omega > 0, \quad t \to \infty.
\end{align}
For these two examples we solve the relevant algebraic system perturbatively up to eighth order. 

The main difference of this new perturbative approach compared with the one used in \cite{trilogy2} is that in the latter approach the analysis was first carried out for all $t$ and then the limit $t \to \infty$ was computed, whereas in the new approach the analysis is carried out directly in the limit of $t \to \infty$. 
The advantage of the approach presented here is that it is computationally much easier and can be used to find coefficients of higher order. The advantage of the approach of \cite{trilogy2} is that it does not require the assumption that the Neumann boundary value is $t$-periodic as $t \to \infty$. However, regarding the latter assumption we note that the new approach gives by construction the large $t$ asymptotics of $g_1(t)$. Thus, if one can use PDE techniques to guarantee existence and uniqueness, the periodicity assumption can be justified a posteriori. 


There exist particular functions $g_0^b(t)$ for which the corresponding functions $g_1^b(t)$ can be obtained in closed form. Among these functions there exists a subclass which have a `linear limit'. For example, among the three pairs defined in (\ref{1.1}) and (\ref{1.3}) the pairs $(\alpha e^{i\omega t}, -\alpha \sqrt{\omega - \alpha^2} e^{i\omega t})$ and $(\alpha e^{i\omega t}, i \alpha \sqrt{|\omega| + 2 \alpha^2} e^{i\omega t})$ have `linear limits'. The above perturbative approach provides a simple and effective way for obtaining $g_1^b(t)$ for this subclass. Indeed, this can be achieved by `summing up' the relevant perturbative series. 
This approach is illustrated using the particular function $g_0^b(t)$ defined by the right-hand side of (\ref{1.1})

For completeness, in addition to the NLS, we also consider the linearized version of the NLS, namely the equation obtained from (\ref{nls}) by deleting the nonlinear term. For this equation, which we call the free Schr\"odinger equation, we derive an expression for the solution in the quarter plane with the same initial and boundary conditions as those specified for the NLS. 
We show that the Neumann boundary value $g_1(t)$ approaches a periodic function as $t \to \infty$. Moreover, we show that this periodic function, which we denote by $g_1^b(t)$, is uniquely determined from $g_0^b(t)$ and we give explicit expressions for its Fourier coefficients. This result is consistent with the perturbative result mentioned earlier for the NLS, since in both cases $g_1^b(t)$ is {\it independent} of $u_0(x)$ in the linear limit.

All results stated in the paper are rigorous. At times, this level of rigor has been achieved at the expense of making assumptions which are stronger than what one would ideally prefer from the point of view of applications. For example, for the NLS we assume that the Dirichlet and Neumann values approach periodic functions at least as fast as $t^{-7/2}$ as $t \to \infty$ (see equation (\ref{uNgN})). This assumption can presumably be weakened. 
For the free Schr\"odinger equation, in order to keep the presentation at a reasonable length, we assume infinite compatibility of the initial and boundary data at the origin. 

\begin{remark}\upshape\label{nonuniqueremark}
For the nonlinear problem, the asymptotic Neumann profile $g_1^b(t)$ is {\it not} uniquely determined from $g_0^b(t)$. This can be seen already from the example (\ref{1.1}) of the single exponential. Indeed, it was shown in \cite{BKS2009} that for the focusing NLS there exists a solution $u$ satisfying $u(0,t) \sim \alpha e^{i\omega t}$ and $u_x(0,t) \sim ce^{i\omega t}$ if and only if the triplet of constants $(\alpha,\omega,c)$ satisfies either
\begin{subequations}\label{1.3}
\begin{align}\label{1.3a}
 c = \pm \alpha\sqrt{\omega - \alpha^2} \quad \text{and} \quad \omega \geq \alpha^2,
\end{align}
or
\begin{align}\label{1.3b}
  c = i \alpha \sqrt{|\omega| + 2 \alpha^2} \quad \text{and}  \quad \omega \leq - 6\alpha^2.
\end{align}
\end{subequations}
Given $\alpha$ and $\omega$, equation (\ref{1.3a}) allows for two different values for $c$. 
Thus $g_1^b(t)$ depends not only on $g_0^b(t)$, but also on $u_0(x)$.
It turns out that among the three branches of solutions given in (\ref{1.3a}) and (\ref{1.3b}), only the branches for which $c = -\alpha\sqrt{\omega - \alpha^2}$ and $c = i \alpha \sqrt{|\omega| + 2 \alpha^2}$ have linear limits. There {\it does} exist a solution of the NLS for which the branch $c = \alpha\sqrt{\omega - \alpha^2}$ occurs (in fact, this branch includes certain stationary soliton solutions), but this solution has {\it no} linear limit, see remark \ref{onesolitonlinearlimitremark}. 
\end{remark}

\section{Linear problem}\nequation\label{linearsec}
Let $\{D_j^0\}_1^4$ denote the four quadrants of the complex $k$-plane:
\begin{align*}
D_1^0 = \{\re k > 0\} \cap \{\im k > 0\}, \qquad D_2^0 = \{\re k < 0\} \cap \{ \im k > 0\}, \\
D_3^0 = \{\re k < 0\} \cap \{\im k < 0\}, \qquad D_4^0 = \{\re k > 0\} \cap \{ \im k < 0\}.
\end{align*}
Let $D_+^0 = D_1^0 \cup D_3^0$ and $D_-^0 = D_2^0 \cup D_4^0$.

\begin{definition}\upshape\label{linearsoldef}
A {\it solution of the free Schr\"odinger equation in the quarter plane} is a smooth function $u:[0,\infty) \times [0,\infty) \to \C$ such that $u(\cdot, t) \in \mathcal{S}([0,\infty))$ for each $t \geq 0$, and such that 
\begin{align}\label{linearnls}
iu_t + u_{xx} = 0
\end{align}
 for $x> 0$ and $t > 0$.
\end{definition}

\begin{proposition}
Suppose $u(x,t)$ is a solution of the free Schr\"odinger equation in the quarter plane and let $u_0(x) = u(x,0)$ and $g_0(t) = u(0,t)$. Then
\begin{align}\nonumber
 u(x,t) = &\; \frac{1}{\pi}\int_\R  e^{-2ikx-f(k) t} \hat{u}_0(k) dk
 	\\ \label{linearsolution}
& + \frac{1}{2i\pi}\int_{\partial D_3^0}  e^{-2ikx-f(k)t} \big[f'(k)\tilde{g}_0(f(k),t)
 + 2i \hat{u}_0(-k)\big]dk,
\end{align}
where, $f(k) = 4i k^2$,
\begin{align} \nonumber
& \hat{u}_0(k) = \int_0^\infty e^{2ikx}u_0(x) dx, \qquad \im k \geq 0,
	\\ \label{linearq0g0def}
& \tilde{g}_0(\kappa,t) =  \int_0^t e^{\kappa s} g_0(s)ds, \qquad \kappa \in \C.
\end{align}
Moreover, the Neumann value $g_1(t) = u_x(0,t)$ is given by
\begin{align}\label{g1linearalt}
g_1(t) = \frac{e^{-\frac{i\pi}{4}}}{\sqrt{\pi}} \bigg[\frac{u_0(0) - g_0(0)}{\sqrt{t}}
+ \frac{1}{\sqrt{t}}\int_0^\infty e^{\frac{ix^2}{4t}}u_{0}'(x) dx
- \int_0^t \frac{\dot{g}_0(s)}{\sqrt{t-s}} ds\bigg], \qquad t > 0,
\end{align}
or, alternatively, by
\begin{align}\nonumber
g_1(t) = &\; \frac{e^{-\frac{i\pi}{4}}}{\sqrt{\pi t}} \bigg(u_0(0) 
+ \int_0^\infty e^{\frac{ix^2}{4t}}u_{0}'(x) dx\bigg)
	\\ \label{g1linear}
& + \frac{2}{\pi} \int_{\partial D_3^0} \bigg(g_0(t) - f(k) \int_0^t e^{f(k)(s-t)}g_0(s)ds\bigg) dk, \qquad t > 0.
\end{align}
\end{proposition}
\proofbegin
Equation (\ref{linearnls}) is the compatibility condition of the Lax pair equations
\begin{align}\label{linearlax}
\begin{cases}
\varphi_x + 2ik \varphi = u,
	\\
\varphi_t + 4i k^2 \varphi = 2ku + i u_x,
\end{cases}
\end{align}
where $k \in \C$ is the spectral parameter and $\varphi(x,t,k)$ is a scalar-valued eigenfunction. 
We write (\ref{linearlax}) in the differential form 
$$d\bigl(e^{2ikx + f(k) t}\varphi\bigr) = W,$$
where the closed one-form $W(x,t,k)$ is defined by
$$W = e^{2ikx + f(k) t}\big[ u dx + (2ku + iu_x)dt\big].$$
Green's Theorem implies that the integral of $W$ around the boundary of the domain $(0,\infty) \times (0,t)$ in the $(x,t)$-plane vanishes. This yields the global relation
\begin{align} \label{linearGR}
& \hat{u}_0(k) - e^{f(k)t}\hat{u}(k,t) + \tilde{g}(k,t) = 0, \qquad  \im k \geq 0,
\end{align}
where 
\begin{align}\nonumber
& \tilde{g}(k,t) = - 2k\tilde{g}_0(f(k),t) - i \tilde{g}_1(f(k),t),
	\\ \nonumber
& 
\hat{u}(k,t) = \int_0^\infty e^{2ikx} u(x,t) dx, \qquad
\tilde{g}_1(\kappa,t) =  \int_0^t e^{\kappa s} g_1(s)ds,
\end{align}
and $\hat{u}_0, \tilde{g}_0$ are defined in (\ref{linearq0g0def}).

Multiplying equation (\ref{linearGR}) by $\frac{1}{\pi}e^{-2ikx-f(k)t}$ and integrating the resulting equation along $\R$ with respect to $k$, we find
\begin{align}\label{solutionformula}
 u(x,t) = &\; \frac{1}{\pi}\int_\R  e^{-2ikx-f(k)t} \hat{u}_0(k) dk
 - \frac{1}{\pi}\int_{\partial D_3^0}  e^{-2ikx-f(k) t} \tilde{g}(k,t)dk,
\end{align}
where we have used Jordan's lemma to deform the contour from $\R$ to $-\partial D_3^0$ in the second integral.

The next step consists of using the global relation to eliminate the unknown function $\tilde{g}_1(k,t)$ from (\ref{solutionformula}). Letting $k \to -k$ in (\ref{linearGR}) and solving for $\tilde{g}_1$, we find 
\begin{align}\nonumber
& \tilde{g}_1(f(k),t) = -i \hat{u}_0(-k) + i e^{f(k)t}\hat{u}(-k,t) - 2ik\tilde{g}_0(f(k),t), \qquad  \im k \leq 0.
\end{align}
Substituting this expression into the solution formula (\ref{solutionformula}) and simplifying, we find (\ref{linearsolution}). Indeed, the contribution from the term involving $\hat{u}(-k,t)$ vanishes since the exponential $e^{-2ikx}$ has decay in $D_3^0$.

In order to derive (\ref{g1linearalt}), we note that (\ref{linearsolution}) implies
\begin{align}\nonumber
 u(x,t) = &\; \frac{1}{\pi}\int_\R  e^{-2ikx-f(k) t} \hat{u}_0(k) dk
+  \frac{1}{\pi}\int_{\partial D_3^0}  e^{-2ikx-f(k) t} \hat{u}_0(-k) dk
 	\\ \nonumber
& + \frac{1}{2i\pi}\int_{\partial \hat{D}_3^0} f'(k) e^{-2ikx-f(k)t} 
\bigg(
\frac{e^{f(k) t}g_0(t) - g_0(0)}{f(k)}  
- \int_0^t \frac{e^{f(k) s}}{f(k)} \dot{g}_0(s)ds\bigg)dk,
\end{align}
where $\partial \hat{D}_3^0$ denotes the contour $\partial D_3^0$ deformed so that it passes to the right of the removable singularity at $k = 0$. 
Since we have deformed the contour to $\partial \hat{D}_3^0$, we can split the $k$-integral and compute the part involving $g_0(t)$ using Cauchy's theorem. We also let $k \to -k$ in the first integral. This gives
\begin{align*}\nonumber
u(x,t) = &\;  \frac{1}{\pi}\int_\R  e^{2ikx-f(k) t} \hat{u}_0(-k) dk
+  \frac{1}{\pi}\int_{\partial D_3^0}  e^{-2ikx-f(k) t} \hat{u}_0(-k) dk
+ 2g_0(t)
	\\
& - \frac{1}{2i\pi}  \int_{\partial \hat{D}_3^0} \frac{f'(k)}{f(k)}  e^{-2ikx - f(k) t}
\bigg(g_0(0) + \int_0^{t} e^{f(k) s} \dot{g}_0(s)ds\bigg) dk.
\end{align*}
Hence, applying $\frac{\partial}{\partial x} \big|_{x = 0}$ to both sides and using Jordan's lemma we find
\begin{align}\nonumber
g_1(t) = &\; \frac{4}{i \pi}\int_{\partial D_3^0} k e^{-f(k) t} \hat{u}_0(-k) dk
	\\ \label{g1linear1}
& + \frac{1}{\pi}  \int_{\partial D_3^0} \frac{k f'(k)}{f(k)}  e^{- f(k) t}
\bigg(g_0(0) + \int_0^{t} e^{f(k) s} \dot{g}_0(s)ds\bigg) dk.
\end{align}
Changing the order of integration in (\ref{g1linear1}) and using the identity
\begin{align}\label{identity}
\int_{\partial D_3^0} e^{-f(k)t - 2ikx} dk = -\frac{e^{-\frac{i\pi}{4}} \sqrt{\pi}}{2\sqrt{t}} e^{\frac{ix^2}{4t}}, \qquad x \geq 0, \quad t > 0,
\end{align}
we find (\ref{g1linearalt}). If we integrate by parts in the second term on the right-hand side of (\ref{g1linear1}) instead, we find (\ref{g1linear}).
\proofend

\begin{remark}\upshape
Definition \ref{linearsoldef} requires that the solution $u$ is smooth $[0,\infty) \times [0,\infty) \to \C$. In particular, this means that $u_0(0) = g_0(0)$. We have kept $u_0(0)$ and $g_0(0)$ separate in the derivation of equations (\ref{g1linearalt}) and (\ref{g1linear}) in order to illustrate the contributions these terms would make in cases when $u_0(0) \neq g_0(0)$. The assumption that $u:[0,\infty) \times [0,\infty) \to \C$ is smooth is made for convenience and is stronger than necessary. Under the assumption of infinite compatibility of the initial and boundary data at the origin, wellposedness for (\ref{linearnls}) in the quarter plane is proved in \cite{FS1999} in the smooth setting. However, it is possible to prove analogous results under weaker assumptions \cite{FHMpreprint}.
\end{remark}

\subsection{Asymptotically admissible pairs for the linear problem}

\begin{definition}\upshape\label{linearadmissibledef}
A pair of smooth functions $\{g_0^b(t), g_1^b(t)\}$, $t \geq 0$, is {\it asymptotically admissible for the free Schr\"odinger equation} if there exists a solution $u(x,t)$ of the free Schr\"odinger equation in the quarter plane (see Definition \ref{linearsoldef}) such that the Dirichlet and Neumann boundary values of $u$ asymptote towards $g_0^b(t)$ and $g_1^b(t)$ respectively in the sense that
\begin{align*}
u(0, \cdot) - g_0^b \in \mathcal{S}([0,\infty)), \qquad u_{x}(0,t) - g_1^b(t) = O(t^{-3/2}), \qquad t\to \infty.
\end{align*}
\end{definition}

\begin{proposition}\label{linearadmissibleprop}
Let 
\begin{align}\label{linearg0B}
  g_0^b(t) = \sum_{n=-\infty}^\infty a_n e^{in\omega t}, \qquad t > 0,
\end{align}
be a smooth periodic function of period $\tau = \frac{2\pi}{\omega} > 0$ with Fourier coefficients $a_n \in \C$, $n \in \Z$, and mean value zero, i.e. $a_0 = 0$.
Then, there exists a unique periodic function  $g_1^b(t)$, such that $\{g_0^b(t), g_1^b(t)\}$ is an asymptotically admissible pair for the free Schr\"odinger equation. This function $g_1^b(t)$ also has period $\tau$ and its Fourier series is given by 
\begin{align}\label{linearg1B}
  g_1^b(t) = \sum_{n=-\infty}^\infty c_n e^{in\omega t}, \qquad t > 0,
\end{align}
where
\begin{align}\label{linearcn}
c_n = \begin{cases} -\sqrt{n\omega} a_n, & n > 0, \\
0, & n = 0, \\
i\sqrt{-n\omega} a_n, & n < 0.
\end{cases}
\end{align}
\end{proposition}
\proofbegin
We first prove that the pair $\{g_0^b, g_1^b\}$ with $g_1^b$ given by (\ref{linearg1B})  is indeed asymptotically admissible. Since $g_0^b$ is smooth, the coefficients $\{a_n\}$ and $\{c_n\}$ converge to zero faster than any power of $n$ as $|n|\to \infty$ and the Fourier series in (\ref{linearg0B}) and (\ref{linearg1B}) converge uniformly in $t$.
Let $u$ be the solution in (\ref{linearsolution}) with Dirichlet datum $g_0(t) = g_0^b(t)$ and any initial datum $u_0 \in \mathcal{S}([0,\infty))$. 
We will show that $g_1(t) = u_x(0,t)$ satisfies $g_1(t) - g_1^b(t) = O(t^{-3/2})$. 
Note that
$$u_0(0) + \int_0^\infty e^{\frac{ix^2}{4t}}u_{0}'(x) dx 
= - \frac{i}{2t} \int_0^\infty x e^{\frac{ix^2}{4t}}u_{0}(x) dx = O(t^{-1}), \qquad t\to\infty.$$
Thus, equation (\ref{g1linear}) implies that the contribution from $u_0(x)$ to $g_1(t)$ is of order $O(t^{-3/2})$.
Noting that $kf'(k) = 2f(k)$, equation (\ref{g1linear}) yields
\begin{align}\label{g1J}
g_1(t) = J(t) + O(t^{-3/2}), \qquad t\to \infty,
\end{align}
where
\begin{align}\nonumber
  J(t) & = \frac{1}{\pi}\int_{\partial D_3^0} kf'(k)\bigg(\frac{g_0(t)}{f(k)} - \int_0^t e^{f(k)(s-t)} g_0(s) ds\bigg) dk
  	\\\nonumber
 & = \frac{1}{\pi}\int_{\partial D_3^0} kf'(k)\sum_{n\neq 0} a_n \bigg(\frac{e^{in\omega t}}{f(k)}  
   - \int_0^t e^{f(k)(s-t)} e^{in\omega s} ds\bigg) dk
   	\\ \label{Jdef}
& =  \frac{1}{\pi}  \int_{\partial D_3^0} kf'(k) \sum_{n\neq 0} a_n \bigg(\frac{e^{in\omega t}}{f(k)}  
   - \frac{e^{in\omega t} - e^{-f(k)t}}{f(k) + in\omega} \bigg) dk.
\end{align}
For each $n \in \Z$, the equation $f(k) + in\omega = 0$ has a unique root $k = K(n)$ in $\partial D_3^0$ given by
$$K(n) = \begin{cases} -\frac{i\sqrt{n\omega}}{2}, & n \geq 0, \\
-\frac{\sqrt{-n\omega}}{2}, & n < 0.
\end{cases}$$
Since the integrand on the right-hand side of (\ref{Jdef}) has removable singularities at the points $\{K(n) \, | \, n \neq 0\}$, we deform the contour to $\partial \hat{D}_3^0$, where $\partial \hat{D}_3^0$ denotes a deformation of $\partial D_3^0$ which passes to the right of these removable singularities. We can then split the integral as follows:
\begin{align}\nonumber
J(t) = &\; \frac{1}{\pi} \sum_{n \neq 0} a_n  e^{in\omega t} \int_{\partial \hat{D}_3^0} kf'(k) \bigg(\frac{1}{f(k)}     - \frac{1}{f(k) + in\omega} \bigg) dk
	\\ \label{J1}
& + 
\frac{1}{\pi}  \int_{\partial \hat{D}_3^0} kf'(k) \bigg( \sum_{n \neq 0} \frac{a_n}{f(k) + in\omega} \bigg)e^{-f(k)t}dk,
\end{align}
where we have used the decay of the $a_n$'s to interchange the order of integration and summation in the first term on the right-hand side of (\ref{J1}). 
Cauchy's theorem implies that the first term on the right-hand side of (\ref{J1}) equals
$$2i \sum_{n \neq 0} a_n e^{in\omega t} \underset{k=K(n)}{\res} k f'(k)  \bigg(\frac{1}{f(k)} - \frac{1}{f(k) + in\omega}\bigg)
= - 2i \sum_{n=-\infty}^\infty a_n K(n) e^{in\omega t} = g_1^b(t).$$
On the other hand, the function $\sum_{n \neq 0} \frac{a_n}{f(k) + in\omega}$ is analytic in a neighborhood of the steepest descent contour $\Gamma = \{u e^{\frac{3i\pi}{4}} | u \in \R\}$. Hence, deforming the contour $\partial \hat{D}_3^0$ to $\Gamma$, a steepest descent argument implies that the second term on the right-hand side of (\ref{J1}) is given by
$$\frac{1}{\pi} \int_{\Gamma} kf'(k) \bigg( \sum_{n\neq 0} \frac{a_n}{f(k) + in\omega} \bigg)e^{-f(k)t}dk
 = O\big(t^{-\frac{3}{2}}\big), \qquad t \to \infty.$$
This estimate together with equations (\ref{g1J}) and (\ref{J1}) yield
$$g_1(t) = g_1^b(t) + O(t^{-3/2}), \qquad t\to \infty.$$
This completes the proof of the admissibility of $\{g_0^b, g_1^b\}$.

In order to prove uniqueness of $g_1^b$, suppose that $\{g_0^b, \tilde{g}_1^b\}$ is another asymptotically admissible pair, where $\tilde{g}_1^b$ is periodic of period $\tilde{\tau} > 0$. Let $u(x,t)$ and $\tilde{u}(x,t)$ be solutions corresponding to $\{g_0^b, g_1^b\}$ and $\{g_0^b, \tilde{g}_1^b\}$ as in Definition \ref{linearadmissibledef}. First note that
\begin{align*}
 \tilde{g}_1^b(t) - g_1^b(t) & = \tilde{g}_1^b(t) - \tilde{u}_x(0,t) + \tilde{u}_x(0,t) - u_x(0,t) + u_x(0,t) - g_1^b(t)
	\\
& = \tilde{u}_x(0,t) - u_x(0,t) + O(t^{-3/2}), \qquad t \to \infty.
\end{align*}
Moreover, the representation (\ref{g1linearalt}) for the Neumann value yields
\begin{align*}
\tilde{u}_x(0,t) - u_x(0,t) = \frac{e^{-\frac{i\pi}{4}}}{\sqrt{\pi}} \int_0^t \frac{\partial_t u(0,s) - \partial_t \tilde{u}(0,s)}{\sqrt{t-s}} ds + O(t^{-\frac{1}{2}})
= O(t^{-\frac{1}{2}}), \qquad t\to \infty,
\end{align*}
where the last step uses the fact that $ u(0,\cdot) - \tilde{u}(0,\cdot) \in \mathcal{S}([0,\infty))$.
Hence $\tilde{g}_1^b(t) - g_1^b(t) = O(t^{-1/2})$. 
This shows that $\tilde{g}_1^b = g_1^b$. Indeed, letting $n \to \infty$ in the following estimate shows that $\tilde{g}_1^b$ has the same period $\tau$ as $g_1^b$:
\begin{align*}
|\tilde{g}_1^b(t_0 + \tau) & - \tilde{g}_1^b(t_0)|
= |\tilde{g}_1^b(t_0 + n \tilde{\tau} +  \tau) - \tilde{g}_1^b(t_0 + n \tilde{\tau})|
	\\
& \leq |\tilde{g}_1^b(t_0 + n \tilde{\tau} +  \tau) - g_1^b(t_0 + n\tilde{\tau} + \tau)| + |g_1^b(t_0 + n\tilde{\tau}) - \tilde{g}_1^b(t_0 + n \tilde{\tau})|
	\\
& = O((t_0 + n \tilde{\tau})^{-1/2}), \qquad n \to \infty, \quad n \in \Z.
\end{align*}
The equation
\begin{align*}
|\tilde{g}_1^b(t_0) - g_1^b(t_0)|
 = |\tilde{g}_1^b(t_0 + n\tau) - g_1^b(t_0 + n \tau)|
\to 0, \qquad n \to \infty, \quad n \in \Z,
\end{align*}
then shows that $\tilde{g}_1^b = g_1^b$.
\proofend

The above uniqueness argument also yields the following result.

\begin{proposition}\label{linearprop}
Suppose $u(x,t)$ is a solution of the free Schr\"odinger equation in the quarter plane (see Definition \ref{linearsoldef}) such that $u(0, \cdot) - g_0^b \in \mathcal{S}([0,\infty))$, where
\begin{align*}
  g_0^b(t) = \sum_{n=-\infty}^\infty a_n e^{in\omega t}, \qquad t > 0,
\end{align*}
is a smooth periodic function of period $\tau = \frac{2\pi}{\omega} > 0$ and mean value zero, i.e. $a_0 = 0$.
Then 
$$u_{x}(0,t) - g_1^b(t) = O(t^{-1/2}), \qquad t\to \infty,$$
where $g_1^b(t)$ is the periodic function given in (\ref{linearg1B})-(\ref{linearcn}). 
\end{proposition}

\begin{remark}[Uniqueness of $g_1^b$ in the linear limit] \upshape
Proposition \ref{linearadmissibleprop} shows that for the free Schr\"odinger equation, the periodic function $g_1^b$ is uniquely determined from $g_0^b$. In particular, the initial datum and the behavior of the Dirichlet datum for small $t$ have no effect on the long-time behavior of the Neumann datum. 
\end{remark}

\begin{remark}\upshape
If the assumption of periodicity on the function $g_1^b(t)$ is dropped in Proposition \ref{linearadmissibleprop}, the uniqueness of $g_1^b(t)$ is clearly lost since any perturbation of $g_1^b(t)$ induced by a small change in $u_0(x)$ gives rise to another asymptotically admissible pair with the same $g_0^b(t)$. However, if we define two pairs $\{g_0^b, g_1^b\}$ and $\{g_0^b, \tilde{g}_1^b\}$ to be equivalent provided that $g_1^b - \tilde{g}_1^b = O(t^{-1/2})$, then uniqueness is valid also without the assumption of periodicity. 
\end{remark}

\subsection{Single exponential}
As an application of the above results, we consider the case where the Dirichlet datum is a single exponential.

\begin{proposition}\label{linearprop2}
Let $u(x,t)$ be any solution of the free Schr\"odinger equation equation (\ref{linearnls}) in the quarter plane such that the Dirichlet datum $g_0(t)$ satisfies $g_0(t) - e^{i\omega t} \in \mathcal{S}([0,\infty))$, $\omega \in \R$. Then the Neumann value $g_1(t) = u_x(0,t)$ satisfies
\begin{align}\label{g1linearasym}
g_1(t) =  O\big(t^{-\frac{1}{2}}\big) + \begin{cases} -\sqrt{\omega} e^{i\omega t}, & \omega > 0, \\
0, & \omega = 0, \\
i\sqrt{|\omega|} e^{i\omega t}, & \omega < 0.
\end{cases}
\end{align}
If $g_0(t) = e^{i\omega t}$, then
\begin{align}\label{g1linearasym2}
g_1(t) =  O\big(t^{-\frac{3}{2}}\big) + \begin{cases} -\sqrt{\omega} e^{i\omega t}, & \omega > 0, \\
-\frac{e^{-\frac{i\pi}{4}}}{\sqrt{\pi t}}, & \omega = 0, \\
i\sqrt{|\omega|} e^{i\omega t}, & \omega < 0.
\end{cases}
\end{align}
\end{proposition}
\proofbegin
The cases $\omega > 0$ and $\omega < 0$ of (\ref{g1linearasym}) follow immediately from proposition \ref{linearprop}. In the case of $\omega = 0$, $g_0^b = 1$ does not have zero mean, hence proposition \ref{linearprop} does not apply. However, the case $\omega = 0$ of (\ref{g1linearasym}) follows from (\ref{g1linearalt}).

The cases $\omega > 0$ and $\omega < 0$ of (\ref{g1linearasym2}) follow from the proof of Proposition \ref{linearadmissibleprop}; the case $\omega = 0$ follows from (\ref{g1linearalt}) and an integration by parts.
\proofend

\section{A new perturbative approach}\nequation \label{pertsec}
The most challenging problem in the analysis of boundary value problems for nonlinear integrable PDEs is the problem of expressing the unknown boundary value(s) in terms of the given initial and boundary conditions. For the Dirichlet problem for the NLS on the half-line, this amounts to determining the Neumann value $g_1(t)$ in terms of the Dirichlet value $g_0(t)$ and the initial datum $u_0(x)$. 
In the context of asymptotically $t$-periodic data, there are situations where we do not necessarily need to know the Neumann value $g_1(t)$ for all $t > 0$, but it is sufficient to know its asymptotic form $g_1^b(t)$. It is therefore natural to ask the following questions: Given the asymptotic form $g_0^b(t)$ of the Dirichlet datum, can we find the asymptotic form $g_1^b(t)$ of the Neumann value? Is $g_1^b(t)$ uniquely determined by $g_0^b(t)$ alone, or does $g_1^b(t)$ also depend on the initial data $u_0(x)$ or on the difference between $g_0(t)$ and $g_0^b(t)$?

The second of these questions is easy to answer: In general, $g_1^b(t)$ is {\it not} uniquely determined from $g_0^b(t)$ alone. Indeed, consider the stationary one-soliton
\begin{align}\label{stationarysoliton}
u(x,t) = \frac{\sqrt{\omega}}{\cosh(x\sqrt{\omega} - \gamma)}  e^{it\omega}, \qquad \gamma \in \R, \quad \omega > 0,
\end{align}
for which
$$g_0(t) = g_0^b(t) = \alpha e^{i\omega t}, \qquad g_1(t) = g_1^b(t) = ce^{i\omega t},$$
with
\begin{align}\label{stationaryc}
\alpha = \frac{\sqrt{\omega}}{\cosh(\gamma)}, 
\qquad c = \sigma \alpha\sqrt{\omega - \alpha^2}, \qquad \sigma = \sgn(\gamma).
\end{align}
Letting $\gamma \to -\gamma$ in (\ref{stationarysoliton}) leaves the Dirichlet datum invariant but switches the sign of the Neumann value. 

However, in the following sections we will show that if a certain boundedness assumption is satisfied, then {\it $g_1^b(t)$ is indeed uniquely determined from $g_0^b(t)$ alone to all orders in a perturbative expansion}. This is consistent with the results of section \ref{linearsec}, where it was shown that $g_1^b(t)$ is uniquely determined from $g_0^b(t)$ in the linear limit.
The stationary one-solitons with $c = \alpha\sqrt{\omega - \alpha^2}$ do not have a linear limit (see Remark \ref{onesolitonlinearlimitremark} below), thus their existence does not contradict the result of this section.

\subsection{Main result}
Before stating the main result, we define what we mean by a perturbative solution.
 
\begin{definition}\upshape
A {\it perturbative solution of the NLS equation (\ref{nls}) in the quarter plane} is a sequence of smooth functions $\{u_N(x,t)\}_1^\infty$ defined for $x>0$ and $t>0$ with the following properties:
\begin{enumerate}[(i)]
\item The formal power series 
\begin{align}\label{uformal}
u = \epsilon u_1 + \epsilon^2 u_2 + \epsilon^3 u_3 + \cdots
\end{align}
satisfies (\ref{nls}) in the quarter plane $\{x> 0, t > 0\}$ to all orders in a perturbative expansion, that is, 
\begin{align*}
O(\epsilon): & \quad iu_{1t} + u_{1xx} = 0,
	\\
O(\epsilon^2): & \quad iu_{2t} + u_{2xx} = 0,
	\\
O(\epsilon^3): & \quad iu_{3t} + u_{3xx} - 2\lambda |u_1|^2 u_1 = 0,
	\\
O(\epsilon^4): & \quad iu_{4t} + u_{4xx} - 2\lambda (u_1^2 \bar{u}_2 + 2u_1 u_2 \bar{u}_1) = 0,
	\\
 \vdots &
	\\
O(\epsilon^N): & \quad iu_{Nt} + u_{Nxx} - 2\lambda \big\{|u|^2 u\big\}_N = 0,
	\\
 \vdots &
\end{align*}
where $\{ \cdots\}_N$ denotes the coefficient of $\epsilon^N$ of the enclosed expression.
	
\item For each $N$, $u_N(x,t)$ and all its partial derivatives have continuous extensions to $\{x\geq 0, t \geq 0\}$.

\item For each $N$ and $t \geq 0$, $u_N(\cdot,t) \in \mathcal{S}([0,\infty))$.
 \end{enumerate}
\end{definition}

The following main result provides a constructive algorithm for computing the asymptotic form of the Neumann value $\{g_{1N}^b(t)\}_{N = 1}^\infty$ from the asymptotic form of the Dirichlet datum $\{g_{0N}^b(t)\}_{N = 1}^\infty$. In other words, it provides an explicit construction of the Dirichlet to Neumann map for asymptotically $t$-periodic data in the limit of large $t$ and small data.

\begin{theorem}\label{pertth}
Let $\{g_{0N}^b(t), g_{1N}^b(t)\}_1^\infty$ be smooth periodic functions of period $\tau = \frac{2\pi}{\omega} > 0$.
Suppose $\{u_N(x,t)\}_1^\infty$ is a perturbative solution of NLS in the quarter plane such that, for each $N$, the Dirichlet and Neumann boundary values of $u_N$ asymptote towards $g_{0N}^b(t)$ and $g_{1N}^b(t)$ respectively in the sense that
\begin{align}\label{uNgN}
\begin{cases} 
u_N(0, t) - g_{0N}^b(t) = O(t^{-7/2}), \\
 u_{Nx}(0,t) - g_{1N}^b(t) = O(t^{-7/2}),
 \end{cases} \qquad t\to \infty, \quad N \geq 1.
\end{align}
Suppose the set $\{u_N(\cdot, t)| t \geq 0\}$ is bounded in $L^1([0,\infty))$ for each $N \geq 1$.\footnote{This assumption can be weakened if $V$ approaches $V^b$ faster than $O(t^{-7/2})$, see Remark \ref{L1remark} below.}

Then the asymptotic Neumann values $\{g_{1N}^b(t)\}_{N = 1}^\infty$ can be constructed explicitly from the asymptotic Dirichlet values $\{g_{0N}^b(t)\}_{N = 1}^\infty$ as follows: Let $\{a_{N,n}\}_{n=-\infty}^\infty$ denote the Fourier coefficients of $g_{0N}^b$:
\begin{align}\label{g0NBFourier}
g_{0N}^b(t) = \sum_{n=-\infty}^\infty a_{N,n} e^{in\omega t}, \qquad N \geq 1.
\end{align}
Then the Fourier coefficients  $\{c_{N,n}\}_{n=-\infty}^\infty$ of
\begin{align}\label{g1NBFourier}
g_{1N}^b(t) = \sum_{n=-\infty}^\infty c_{N,n} e^{in\omega t}, \qquad N \geq 1,
\end{align}
are given by
\begin{align}\nonumber
c_{N,n} = &\; \bigg\{ \bigg(
-2i\lambda k \sum_{l,m=-\infty}^\infty \bar{a}_l d_m(k) d_{n+l-m}(k) 
 - \lambda  \sum_{l, m=-\infty}^\infty  \bar{c}_{l} d_m(k) d_{n+l-m}(k)
	\\ \label{cNndef}
&+ 2\lambda \sum_{l,m = -\infty}^\infty \bar{a}_l a_m d_{n+l-m}(k) + 2ik a_n\bigg)\bigg|_{k = k_1(n)} \biggr\}_N,
\end{align}
where 
\begin{itemize}
\item $a_n = \sum_{N=1}^\infty a_{N, n} \epsilon^N$, $c_n = \sum_{N=1}^\infty c_{N, n} \epsilon^N$, and 
$d_{n}(k) = \sum_{N=1}^\infty d_{N, n}(k) \epsilon^N$. 

\item The coefficients $d_{N,n}(k)$, $N \geq 1$, $n \in \Z$, satisfy
\begin{align}\nonumber
& d_{N,n}(k) = \frac{1}{4ik^2 + in\omega}\biggl\{
-2\lambda k \sum_{l,m=-\infty}^\infty \bar{a}_l d_m(k) d_{n+l-m}(k) 
	\\ \nonumber
& + \lambda i \sum_{l, m=-\infty}^\infty  \bar{c}_{l} d_m(k) d_{n+l-m}(k)
- 2i\lambda \sum_{l,m = -\infty}^\infty \bar{a}_l a_m d_{n+l-m}(k) + 2k a_n + i c_n\biggr\}_N, 
	\\ \label{dNndef}
& \hspace{7cm} n \in \Z, \quad N \geq 1, \quad k \in \partial D_1^0.
\end{align}

\item $k_1(n)$ denotes the unique root of $4k^2 + n\omega = 0$ in $\partial D_1^0$, i.e.
$$k_1(n) = \begin{cases} \frac{i\sqrt{n\omega}}{2}, & n \geq 0, \\
\frac{\sqrt{-n\omega}}{2}, & n < 0.
\end{cases}$$
\end{itemize}
\end{theorem}

The proof of theorem \ref{pertth} will be presented in section \ref{proofsec}. In the remainder of this section, we explain how the $c_{N,n}$'s are determined from equations (\ref{cNndef})-(\ref{dNndef}) and we also provide several examples. 

\subsection{Construction of the $c_{N,n}$'s}
Equation (\ref{cNndef}) with $N = 1$ yields
\begin{subequations}\label{c1nd1n}
\begin{align}\label{c1n}
  c_{1,n} = 2i k_1(n) a_{1,n}, \qquad n \in \Z.
\end{align}
Substituting this into equation (\ref{dNndef}) with $N = 1$, we find
\begin{align}\label{d1n}
d_{1,n}(k) = \frac{2k a_{1,n} + i c_{1,n}}{4ik^2 + in\omega} 
= \frac{a_{1,n}}{2i(k+ k_1(n))}, 
\qquad n \in \Z, \quad k \in \partial D_1^0.
\end{align}
\end{subequations}
Similarly, equations (\ref{cNndef}) and (\ref{dNndef}) with $N = 2$ yield
\begin{subequations}\label{c2nd2n}
\begin{align}\label{c2n}
& c_{2,n} = 2i k_1(n) a_{2,n}, && n \in \Z,
	\\ \label{d2n}
&d_{2,n}(k)  = \frac{a_{2,n}}{2i(k+ k_1(n))}, && n \in \Z, \quad k \in \partial D_1^0.
\end{align}
\end{subequations}
Continuing in this way, equations (\ref{cNndef}) and (\ref{dNndef}) with $N = 3$ yield
\begin{subequations}\label{c3nd3n}
\begin{align} \nonumber
c_{3,n} = &\; \bigg(
-2i\lambda k \sum_{l,m=-\infty}^\infty \bar{a}_{1,l} d_{1,m}(k) d_{1,n+l-m}(k) 
 - \lambda \sum_{l, m=-\infty}^\infty  \bar{c}_{1,l} d_{1,m}(k) d_{1,n+l-m}(k)
	\\ \label{c3n}
&+ 2\lambda \sum_{l,m = -\infty}^\infty \bar{a}_{1,l} a_{1,m} d_{1,n+l-m}(k) + 2i k a_{3,n} \bigg)\bigg|_{k = k_1(n)}, \qquad n \in \Z,
	\\ \nonumber
d_{3,n}(k) = &\; \frac{1}{4ik^2 + in\omega}\biggl(
- 2\lambda k \sum_{l,m=-\infty}^\infty \bar{a}_{1,l} d_{1,m}(k) d_{1,n+l-m}(k) 
	\\ 
&+ \lambda i \sum_{l, m=-\infty}^\infty  \bar{c}_{1,l} d_{1,m}(k) d_{1,n+l-m}(k)
	\\ \nonumber
& - 2i\lambda \sum_{l,m=-\infty}^\infty \bar{a}_{1,l} a_{1,m} d_{1,n+l-m}(k)
 + 2 k a_{3,n} +i c_{3,n}\biggr), \qquad n \in \Z, \quad k \in \partial D_1^0.
\end{align}
\end{subequations}
This process can be continued indefinitely. Indeed, suppose we have determined $\{c_{M,n}\}_{n=-\infty}^\infty$ and $\{d_{M,n}(k)\}_{n=-\infty}^\infty$ for $N \leq M-1$. Then equation (\ref{cNndef}) with $N = M$ yields
\begin{align*}
c_{M,n} = \bigg(F_{Mn}(k) + 2 ik a_{M,n}\bigg)\bigg|_{k = k_1(n)}, \qquad n \in \Z,
\end{align*}
where the function $F_{Mn}(k)$ is given in terms of known lower order terms:
\begin{align*}
F_{Mn}(k) = &\; \lambda  \sum_{l,m=-\infty}^\infty \bigg\{-2i k \bar{a}_l d_m(k) d_{n+l-m}(k) 
 -   \bar{c}_{l} d_m(k) d_{n+l-m}(k)
	\\
&+ 2 \bar{a}_l a_m d_{n+l-m}(k)\bigg\}_M.
\end{align*}
We can now use equation (\ref{dNndef}) with $N = M$ to determine $d_{M,n}(k)$:
\begin{align*}
d_{M,n}(k) = &\;\frac{1}{4ik^2 + in\omega}\bigl(
-iF_{Mn}(k) + 2 k a_{M,n} + i c_{M,n} \bigr), \qquad n \in \Z.
\end{align*}
This determines $\{c_{M,n}\}_{n=-\infty}^\infty$ and $\{d_{M,n}(k)\}_{n=-\infty}^\infty$ for $N = M$ and completes the inductive step.

\subsection{Examples}

\begin{example}[Single exponential]\upshape\label{singleexample1}
Suppose 
$$g_0^b(t) = \epsilon e^{i\omega t}, \qquad t \geq 0, \quad \omega > 0.$$
In this case, all coefficients $a_{N,n}$ are zero except for $a_{1,1} = 1$.
Equations (\ref{c1nd1n}) imply that all the coefficients $c_{1,n}$ and $d_{1,n}(k)$ vanish except for
\begin{align*}
  c_{1,1} = - \sqrt{\omega},
\qquad
d_{1,1}(k) = -\frac{i}{2k + i\sqrt{\omega}}.
\end{align*}
Equations (\ref{c2nd2n}) yield $c_{2,n} = d_{2,n}(k) = 0$ for all $n$.
Equations (\ref{c3nd3n}) imply that all the coefficients $c_{3,n}$ and $d_{3,n}(k)$ vanish except for
\begin{align} \nonumber
c_{3,1} = - \frac{\lambda}{2\sqrt{\omega}},
\qquad
d_{3,1}(k) = \frac{\lambda}{2(2ik - \sqrt{\omega})^2 \sqrt{\omega}}.
\end{align}
Continuing in this way, we find that the nonzero coefficients $c_{N,n}$ with $N \leq 8$ are
$$c_{1,1} = - \sqrt{\omega}, \qquad
c_{3,1} = - \frac{\lambda}{2\sqrt{\omega}}, \qquad
c_{5,1} = \frac{1}{8\omega^{3/2}}, \qquad
c_{7,1} = - \frac{\lambda}{16 \omega^{5/2}}.$$
In summary, we have found that
$$g_1^b(t) = -\epsilon \biggl( \sqrt{\omega} + \frac{\epsilon^2 \lambda }{2\sqrt{\omega}} - \frac{\epsilon^4}{8\omega^{3/2}} + \frac{\epsilon^6 \lambda}{16 \omega^{5/2}} + O(\epsilon^8)\biggr)e^{i\omega t}.$$
The summation of this perturbative expansion suggests
\begin{align}\label{g1Bsingleexp}
g_1^b(t) = -\epsilon \sqrt{\omega + \lambda \epsilon^2} e^{i\omega t},
\end{align}
which, upon identifying $\alpha$ and $\epsilon$, is in agreement with (\ref{1.3a}) (note that the inequality $\omega \geq \alpha^2$ is automatically satisfied in the perturbative limit $\alpha \to 0$). 
\end{example}

\begin{remark}\upshape\label{onesolitonlinearlimitremark}
In addition to the branch $c = - \alpha\sqrt{\omega - \alpha^2}$ found in (\ref{g1Bsingleexp}), equation (\ref{1.3a}) also allows for $c = \alpha\sqrt{\omega - \alpha^2}$. The reason the latter branch of solutions does not show up in example \ref{singleexample1} is that the associated solutions of NLS do not have a linear limit which decays as $x \to \infty$. This can be understood in terms of the stationary one-solitons. Indeed, according to (\ref{stationaryc}) the triple $(\alpha,\omega, c)$ associated with the soliton (\ref{stationarysoliton}) satisfies $c = \sigma \alpha\sqrt{\omega - \alpha^2}$ where $\sigma = \sgn \gamma$ (we assume that $\gamma \neq 0$; if $\gamma = 0$, then $\omega = \alpha^2$ and hence $c = 0$). 
Writing the solution (\ref{stationarysoliton}) in terms of $\alpha = \frac{\sqrt{\omega}}{\cosh(\gamma)} > 0$ as
$$u(x,t) = \frac{\alpha \sqrt{\omega } e^{i\omega t}}{\sqrt{\omega } \cosh \left(x \sqrt{\omega}\right) - \sigma \sqrt{\omega -\alpha^2} \sinh \left(x \sqrt{\omega }\right)},$$	
we find the following expansion as $\alpha \to 0$:
\begin{align*}
u(x,t) = e^{\sigma x \sqrt{\omega } + i\omega t} \bigg(\alpha + \alpha^3 \frac{1 - e^{\sigma 2 x \sqrt{\omega }}}{4 \omega }
+ \alpha^5 \frac{4 -6 e^{\sigma 2 x \sqrt{\omega }} + 2 e^{\sigma 4 x \sqrt{\omega }}}{32 \omega ^2}  + O(\alpha^7)\bigg).
   \end{align*}
The coefficients of this expansion decay as $x \to \infty$ if $\sigma = - 1$, but {\it not} if $\sigma = 1$. 
Consequently, the coefficients constitute a perturbative solution $\{u_N(x,t)\}_1^\infty$ of the NLS on the half-line if $\gamma < 0$, but not if $\gamma > 0$. In other words, the one-solitons associated with $c = - \alpha\sqrt{\omega - \alpha^2}$ have perturbative limits, whereas those associated with $c = \alpha\sqrt{\omega - \alpha^2}$ do not. 
\end{remark}

\begin{example}[Single exponential]\upshape
Suppose 
$$g_{01}^b(t) = \epsilon e^{-i\omega t}, \qquad t \geq 0, \quad \omega > 0.$$
In this case, all coefficients $a_{N,n}$ are zero except for $a_{1,-1} = 1$.
Proceeding as in the previous example, we find that the nonzero coefficients $c_{N,n}$ with $N \leq 8$ are
$$c_{1,-1} = i \sqrt{\omega}, \qquad
c_{3,-1} = - \frac{i \lambda}{\sqrt{\omega}}, \qquad
c_{5,-1} = - \frac{i}{2\omega^{3/2}}, \qquad
c_{7,-1} = - \frac{i \lambda}{2 \omega^{5/2}}.$$
In summary, 
$$g_1^b(t) = i \epsilon \biggl(\sqrt{\omega} - \frac{\epsilon^2 \lambda }{\sqrt{\omega}} - \frac{\epsilon^4}{2\omega^{3/2}} - \frac{\epsilon^6 \lambda}{2 \omega^{5/2}} + O(\epsilon^8)\biggr)e^{-i\omega t}.$$
The summation of this perturbative expansion suggests
$$g_1^b(t) = i \epsilon \sqrt{\omega - 2 \lambda \epsilon^2} e^{-i\omega t},$$
which, upon identifying $\alpha$ and $\epsilon$ and letting $\omega \to -\omega$, is in agreement with (\ref{1.3b}) (note that the inequality $\omega \leq -6\alpha^2$ is automatically satisfied in the perturbative limit $\alpha \to 0$). 
\end{example}

\begin{example}[Sum of exponentials]\upshape\label{sumexample}
We consider the case of
\begin{align}\label{Dirichletcondition}
g_0^b(t) = \epsilon (\alpha e^{i\omega t} + \beta e^{-i\omega t}), \qquad t \geq 0,
\end{align}
where $\alpha, \beta \in \C$ and $\omega > 0$ are constants. In this case, all coefficients $a_{N,n}$ vanish except for 
$$a_{1,1} =  \alpha, \qquad a_{1,-1} = \beta.$$
Equations (\ref{c1nd1n}) imply that all the $c_{1,n}$'s and $d_{1,n}(k)$'s vanish except for
\begin{align*}
c_{1,1} = -\alpha \sqrt{\omega}, \quad 
c_{1,-1} = i \beta \sqrt{\omega}, \quad
d_{1,1}(k) = -\frac{i\alpha}{2k + i\sqrt{\omega}}, \quad
d_{1,-1}(k) = -\frac{i\beta}{2k + \sqrt{\omega}}.
\end{align*}
Equations (\ref{c2nd2n}) yield $c_{2,n} = d_{2,n}(k) = 0$ for all $n$.
In general, $c_{N,n} = 0$ unless both $N$ and $n$ are odd and $N \geq n$.
Continuing in this way, we find that the nonzero coefficients $c_{N,n}$ with $N \leq 8$ are
\begin{subequations}\label{cmnexpressions}
\begin{align}\label{c1expressions}
  c_{1,1} = - \alpha \sqrt{\omega}, \qquad c_{1,-1} = i\beta \sqrt{\omega},
\end{align}  
\begin{align} \nonumber
 & c_{3,3} = -\frac{i (\sqrt{3}+(-2-i)) \alpha ^2 \bar{\beta} \lambda }{2
   \sqrt{\omega }}, &&
  c_{3,1} = -\frac{\alpha  \lambda  (|\alpha|^2 +4 |\beta|^2)}{2 \sqrt{\omega }},
  	\\ \label{c3expressions}
&  c_{3,-1} = -\frac{i \beta  \lambda  (|\beta|^2 +(1-i) |\alpha|^2)}{\sqrt{\omega }}, &&
  c_{3,-3} = \frac{(\sqrt{3}-2-i) \beta ^2 \bar{\alpha} \lambda }{2
   \sqrt{\omega }},
\end{align}
\begin{align} \nonumber
&  c_{5,5} = \frac{((2+i)-(5+2 i) \sqrt{3}-(2+3 i) \sqrt{5}+(3+2 i) \sqrt{15})
   \alpha ^3 \bar{\beta}^2}{16 \omega ^{3/2}}, 
   	\\ \nonumber
&  c_{5,3} = -\frac{i \alpha^2 \bar{\beta} (((3+6 i)+(2+3 i)
   \sqrt{3}) |\alpha|^2+(5+2 i) ((3+2 i)+(2+i)
   \sqrt{3}) |\beta|^2)}{2 ((12+3 i)+(7+2 i)
   \sqrt{3}) \omega ^{3/2}} , 
  	\\\nonumber
&  c_{5,1} = \frac{\alpha  (|\alpha|^4 -4 (\sqrt{3}-6)
   |\alpha|^2  |\beta|^2+2 (9-i   \sqrt{3}) |\beta|^4)}{8 \omega ^{3/2}}, 
   	\\\nonumber
&  c_{5,-1} = \frac{\beta ((\sqrt{3}+(-2+3 i)) |\alpha|^4 - 2 (\sqrt{3}+(2-5 i)) |\alpha|^2  |\beta|^2 -2 i |\beta|^4)}{4 \omega ^{3/2}}, 
  	\\\nonumber
&  c_{5,-3} = -\frac{\beta^2 \bar{\alpha} (((48+9 i)+(29+4 i) \sqrt{3}) |\alpha|^2+8 ((3+5 i)+(2+3 i)   \sqrt{3}) |\beta|^2)}{4 ((21+12 i)+(12+7 i) \sqrt{3}) \omega ^{3/2}} , 
   	\\ \label{c5expressions}
& c_{5,-5} = \frac{((1-2 i)-(2-5 i) \sqrt{3}-(3-2 i) \sqrt{5}+(2-3 i) \sqrt{15})
   \beta ^3 \bar{\alpha}^2}{16 \omega ^{3/2}},
\end{align}
\begin{align} \nonumber 
c_{7,7} = &\; \frac{\alpha ^4 \bar{\beta}^3 \lambda}{32 (56+23
   \sqrt{7}) \omega ^{5/2}}
   \Big\{(-595+147 i)+(525-147 i) \sqrt{3}+(196+49 i) \sqrt{5}
   	\\ \nonumber
&   -(223-30 i) \sqrt{7}-(189-14 i) \sqrt{15}+(192-39 i) \sqrt{21}+(67+19 i) \sqrt{35}
	\\ \nonumber
& -(63+i)  \sqrt{105}\Big\}, 
   	\\ \nonumber
c_{7,5} = &\; \frac{i \alpha ^3 \bar{\beta}^2 \lambda}{8 ((26+97 i)+(15+56 i) \sqrt{3})
   ((25+10 i)+(11+4 i) \sqrt{5}) \omega ^{5/2}}
   	\\ \nonumber
&\times     \Big\{\big((1110+2345 i)+(645+1350 i) \sqrt{3}+(106+1011 i) \sqrt{5}
	\\ \nonumber
& +(63+582 i) \sqrt{15}\big) |\alpha|^2
+(1+i) \big((4542+903 i)+(2623+524 i) \sqrt{3}
	\\ \nonumber
& +(1612+1359 i) \sqrt{5}+(931+786 i) \sqrt{15}\big) |\beta|^2 \Big\}, 
   	\\ \nonumber
 c_{7,3} =  & -\frac{\alpha ^2 \bar{\beta} \lambda}{96 ((627+2340 i)+(362+1351 i)
   \sqrt{3}) \omega ^{5/2}}   \Big\{\big((42588+1476 i)
   	\\ \nonumber
&+(24588+852 i) \sqrt{3}\big) |\alpha|^4 
-(6-6 i) \big((-25328-35261 i)
	\\ \nonumber
& -(14623+20358 i) \sqrt{3}+(3594+1545 i) \sqrt{5}+(2075+892 i) \sqrt{15}\big) |\alpha|^2  |\beta|^2
	\\ \nonumber
& +\big((395885+21770 i)+(228564+12569 i)
   \sqrt{3}-(9405+7020 i) \sqrt{5}
   	\\ \nonumber
&   -(5430+4053 i) \sqrt{15}\big) |\beta|^4 - 96 i \big((795+627 i)+(459+362 i) \sqrt{3}\big)
   |\alpha|^2  |\beta|^2\Big\}, 
  	\\\nonumber
c_{7,1} = & -\frac{\alpha  \lambda }{48 \omega ^{5/2}} 
\Big\{3 |\alpha|^6+6 (4+5 \sqrt{3}) |\alpha|^4  |\beta|^2 +2 \big((162+51 i)
	\\ \nonumber
&-(13+12 i) \sqrt{3}\big) |\alpha|^2 |\beta|^4 
 - 6  |\alpha|^2  |\beta|^2\big((\sqrt{3}-6) |\alpha|^2
   	\\ \nonumber
& +(2+2 i)  (\sqrt{3}+(2+5 i)\big) |\beta|^2)  +16 (15-2 i \sqrt{3}) |\beta|^6\Big\}, 
   	\\\nonumber
c_{7,-1} = &\; \frac{\beta  \lambda}{24 \omega ^{5/2}} \Big\{3 (3-2 i \sqrt{3})|\alpha|^6 +\big((39+42 i)+(28-45 i) \sqrt{3}\big) |\alpha|^4|\beta|^2
	\\ \nonumber
& +((36-15 i)+(23+36 i) \sqrt{3}) |\alpha|^2|\beta|^4 +3
   |\alpha|^2 |\beta|^2 \big((-2+2 i) (\sqrt{3}-6)  |\alpha|^2
   	\\ \nonumber
&   -(\sqrt{3}+(2+5 i)) |\beta|^2 \big)-12 i |\beta|^6\Big\}, 
  	\\\nonumber
c_{7,-3} = &\; \frac{(1+i) \beta^2 \bar{\alpha} \lambda }{96((2340+8733 i)+(1351+5042 i) \sqrt{3}) \omega ^{5/2}} \Big\{3 \big((32088+120795 i)
   	\\ \nonumber
& +(18526+69741 i) \sqrt{3}  +(8106-6393  i) \sqrt{5}+(4680-3691 i) \sqrt{15}\big) |\alpha|^4
   	\\ \nonumber
& +6  i ((168629+32048 i)+(97358+18503 i) \sqrt{3}
 +(8733+11700 i) \sqrt{5}
 	\\ \nonumber
& +(5042+6755 i) \sqrt{15})|\alpha|^2 |\beta|^2 
-(4-4 i) \beta \big(4 ((265+7718 i)
	\\ \nonumber
& +(153+4456 i) \sqrt{3}) |\beta|^2 \bar{\beta}-9 ((892+3329 i)+(515+1922 i)  \sqrt{3}) |\alpha|^2 \bar{\beta}\big)\Big\}, 
   	\\ \nonumber
c_{7,-5} = & -\frac{(1-i) \beta^3 \bar{\alpha}^2 \lambda  }{24((97+26 i)+(56+15 i) \sqrt{3}) ((15+40 i)+(7+18 i) \sqrt{5}) \omega ^{5/2}}
   	\\ \nonumber
&\times  \Big\{3 \big((3820+3375 i)+(2205+1950 i) \sqrt{3}+(1464+1395 i)
   \sqrt{5}
   	\\ \nonumber
&   +(845+806 i) \sqrt{15}\big) |\alpha|^2+(1+i) \big((7755+7962 i)+(4466+4602 i) \sqrt{3}
   	\\ \nonumber
&   +(2742+2586 i) \sqrt{5}+(1578+1495 i) \sqrt{15}\big) |\beta|^2\Big\}, 
   	\\ \nonumber
c_{7,-7} = &\; \frac{\beta ^4 \bar{\alpha}^3 \lambda}{32 (56+23 \sqrt{7}) \omega ^{5/2}}
   \Big\{(-147-595 i)+(147+525 i) \sqrt{3}-(49-196 i) \sqrt{5}
   	\\ \nonumber
&   -(30+223 i) \sqrt{7}-(14+189 i) \sqrt{15}+(39+192 i) \sqrt{21}-(19-67 i) \sqrt{35}
	\\ \label{c7expressions}
& +(1-63  i) \sqrt{105}\Big\} .	
\end{align}
\end{subequations}
\end{example}

\begin{example}[Sine wave] \upshape
It was proved in Theorem 3.2 of \cite{trilogy2} that if $g_0(t) = \epsilon g_{01}(t)$ with $g_{01}(t) = \sin t$, then
\begin{align} \label{g13trilogy2}
& \lambda g_{13}(t) = \mathcal{A} \cos{t} + \mathcal{B} \sin{t} + \frac{(1-\sqrt{3}) (\sqrt{3}+i)}{16}(\cos{3t} - \sin{3t}) + O(t^{-1/16}), && t\to \infty,
\end{align}
where $\mathcal{A},\mathcal{B} \in \C$ are constants. In \cite{trilogy2} rather complicated integral expressions for $\mathcal{A},\mathcal{B}$ were derived. Using the results of this section, we can now determine $\mathcal{A},\mathcal{B}$ explicitly. Indeed, letting $\alpha = \frac{1}{2i}$ and $\beta = -\frac{1}{2i}$ in (\ref{c3expressions}), we find
\begin{align*}
  \mathcal{A} = \frac{1}{4}+\frac{3 i}{16}, \qquad
  \mathcal{B} = -\frac{7}{16}-\frac{i}{4}.
\end{align*}
\end{example}

\begin{example}[Stationary solitons]\upshape
The focusing NLS equation admits the following family of $N$-soliton solutions (see \cite{PZ2004}, Section 1.7.1):
$$u(x,t) = \frac{\det R(x,t)}{\det M(x,t)},$$
where the $N\times N$ matrix $M(x,t)$ is defined by
$$M_{n,k} = \frac{1 + \bar{g}_kg_n}{\bar{\lambda}_n - \lambda_k}, \qquad n,k = 1, \dots, N,$$
the $(N+1) \times (N+1)$ matrix $R(x,t)$ is given by
\begin{align*}
&R_{n,k} = M_{n,k},\qquad n,k = 1, \dots, N; && R_{n,N+1} = g_n, \qquad n = 1, \dots, N;
	\\
&R_{N+1,k} = 1, \qquad k = 1, \dots, N; && R_{N+1,N+1} = 0,
\end{align*}
the scalar-valued function $g_n(x,t)$ is defined by $g_n(x,t) = \gamma_n e^{i(\lambda_n x - \lambda_n^2 t)}$, and $\{\lambda_n\}_1^N$ and $\{\gamma_n\}_1^N$ are arbitrary complex constants such that $\im \lambda_n > 0$, $\lambda_n \neq \lambda_k$ for $n \neq k$, and $\gamma_n \neq 0$.
If the $\lambda_n$'s belong to $i[0,\infty)$ and the quotients $\lambda_n/\lambda_m$ for $1 \leq n,m \leq N$ are rational numbers, then the solution is $t$-periodic. 

Consider the particular case of
$$N = 2, \quad \lambda_1 = i, \quad \lambda_2 = 2i, \quad \gamma_1 = \gamma_2= -i\epsilon,$$
so that
$$u(x,t) = -\frac{6 \epsilon e^{x+4 i t} \left(-2 \epsilon^2 e^{x+3 i t}+\epsilon^2-2 e^{3 x+3 i t}+e^{4
   x}\right)}{\epsilon^4 e^{3 i t}+\epsilon^2 e^{2 x} \left(-8 e^{x+6 i t}+9 e^{2 x+3 i t}+9
   e^{3 i t}-8 e^x\right)+e^{6 x+3 i t}}.$$
We will use the results of theorem \ref{pertth} to reconstruct the Neumann value of this solution from its Dirichlet datum
\begin{align*}
g_0(t) = &\; \frac{6 \epsilon \left(\epsilon^2+1\right) e^{4 i t} \left(-1+2 e^{3 i t}\right)}{\epsilon^4 e^{3 i
   t}-2 \epsilon^2 \left(-9 e^{3 i t}+4 e^{6 i t}+4\right)+e^{3 i t}}
   	\\
= &\; (12 e^{4 i t}- 6e^{i t}) \epsilon + (198 e^{i t}- 48 e^{-2 i t}- 252 e^{4 i t}+ 96 e^{7 i t}) \epsilon^3 + O(\epsilon^5), \qquad \epsilon \to 0.
\end{align*}
Since both $u(x,t)$ and $u_x(x,t)$ are $t$-periodic with period $2\pi$, we have $g_0^b = g_0$ and $g_1^b = g_1$. 
The nonzero coefficients $a_{N,n}$ with $N \leq 4$ are
\begin{align*}
  a_{1,1} = -6, \quad a_{1,4} = 12, \quad
  a_{3,-2} = -48, \quad a_{3,1} = 198, \quad
  a_{3,4} = -252, \quad a_{3,7} = 96.
\end{align*}
Equations (\ref{c1nd1n}) imply that all the $c_{1,n}$'s and $d_{1,n}(k)$'s vanish except for
\begin{align*}
c_{1,1} = 6, \qquad
c_{1,4} = - 24, \qquad
d_{1,1}(k) = \frac{6i}{2k + i}, \qquad
d_{1,4}(k) = - \frac{6i}{k + i}.
\end{align*}
Equations (\ref{c2nd2n}) yield $c_{2,n} = d_{2,n}(k) = 0$ for all $n$.
Equation (\ref{c3n}) implies that all coefficients $c_{3,n}$ vanish except for
\begin{align*}
c_{3,-2} = 192, \qquad c_{3,1} = -882, \qquad
c_{3,4} = 1224, \qquad c_{3,7} = -480 .
\end{align*}
This yields
\begin{align*}
g_1(t) = &\; ( 6e^{i t}- 24 e^{4 i t}) \epsilon + (192 e^{-2 i t}-882 e^{i t}+1224 e^{4 i t}-480 e^{7 i t}) \epsilon^3 + O(\epsilon^5), \quad \epsilon \to 0,
\end{align*}
which is indeed the correct expansion of the Neumann value.
\end{example}

\section{Eigenfunctions}\nequation
In preparation for the proof of Theorem \ref{pertth}, we introduce three eigenfunction solutions $\{\phi_j(x,t,k)\}_1^3$ of the Lax pair of (\ref{nls}) which is given by
\begin{align}\label{lax}
\begin{cases}
  \phi_x + ik\sigma_3 \phi = U\phi,
  	\\ 
  \phi_t + 2ik^2 \sigma_3 \phi = V \phi,	
\end{cases}
\end{align}
where
\begin{align*}
U = \begin{pmatrix} 0 & u \\
\lambda \bar{u} & 0 \end{pmatrix}, \qquad
V = \begin{pmatrix} -i\lambda |u|^2 & 2ku + iu_x \\
2\lambda k \bar{u} - i\lambda \bar{u}_x & i\lambda |u|^2 \end{pmatrix}, \qquad
\sigma_3 = \begin{pmatrix} 1 & 0 \\ 0 & -1 \end{pmatrix}.
\end{align*}

The solutions $\phi_2(x,t,k)$ and $\phi_3(x,t,k)$ of (\ref{lax}) are normalized at $(x_2, t_2) = (0,0)$, and $(x_3, t_3) = (\infty, t)$ respectively. Namely, these eigenfunctions are defined by
$$\phi_j(x,t,k) = \mu_j(x,t,k) e^{-i(kx + 2k^2 t)\sigma_3}, \qquad j = 2,3,$$
where $\{\mu_j(x,t,k)\}_2^3$ are the unique solutions of the linear Volterra integral equation
\begin{align}\label{mu23def}
  \mu_j(x,t,k) = &\; I + \int_{(x_j, t_j)}^{(x,t)} e^{i[k(x'-x) + 2k^2(t'-t)] \hat{\sigma}_3} W_j(x',t',k), \qquad j = 2,3,
\end{align}
with 
$$W_j = (U dx + V dt)\mu_j.$$

The solution $\phi_1(x,t,k)$ of (\ref{lax}) is defined as follows (we refer to \cite{tperiodicI} for further details).
Let $V^b$ denote the function obtained by replacing $u(0,t)$ and $u_x(0,t)$ by $g_0^b(t)$ and $g_1^b(t)$, respectively, in the expression for $V(0,t,k)$.
Let $\psi(t,k)$ be the solution of the background $t$-part
\begin{align}\label{tpartB}
  \psi_t + 2ik^2\sigma_3 \psi = V^b\psi,
\end{align}  
normalized by $\psi(0,k) = I$. We define the entire $2\times 2$-matrix valued function $Z(k)$ by $Z(k) = \psi(\tau,k)$.
The eigenvalues of $Z(k)$ are given by $z(k)$ and $z(k)^{-1}$ where 
\begin{align}\label{lambdaGdef}
z(k) = \frac{1}{2} \big(\tr Z(k) - \sqrt{G(k)}\big), \qquad G(k) = (\tr Z(k))^2 - 4.
\end{align}
Let $\mathcal{P}$ denote the set of branch points defined by
\begin{align}\label{branchpts}
\mathcal{P} = \{k \in \C \; | \; G(k) = 0,  \; \text{or} \; Z_{12}(k) = 0, \; \text{or} \; Z_{21}(k) = 0\},
\end{align}
where $Z_{ij}$ denotes the $(ij)$'th entry of $Z$. 
The set $\mathcal{P}$ is the union of a finite number of zero sets of entire functions, thus $\mathcal{P}$ is a countable set without accumulation points. 

Let $\mathcal{C}$ denote a set of branch cuts connecting all points in $\mathcal{P}$. We choose these branch cuts so that $\mathcal{C}$ is invariant under the involution $k \mapsto \bar{k}$. Letting
\begin{align}\label{mathcalSdef}
S^b(k) = \sqrt{-\frac{Z_{11} - Z_{22} - \sqrt{G}}{2\sqrt{G}}} 
\begin{pmatrix} 1 & - \frac{2 Z_{12}}{Z_{11}-Z_{22}-\sqrt{G}} \\
 \frac{2 Z_{21}}{Z_{11}-Z_{22}-\sqrt{G}} & 1 \end{pmatrix}, \qquad k \in \C \setminus \mathcal{C},
 \end{align}
we find that $S^b(k)$ has unit determinant and that
$$Z(k) = S^b(k) \begin{pmatrix} z(k) & 0 \\ 0 & z^{-1}(k) \end{pmatrix} S^b(k)^{-1}, \qquad k \in \C \setminus \mathcal{C}.$$
The identity
$$(Z_{11} - Z_{22} - \sqrt{G})(Z_{11} - Z_{22} + \sqrt{G}) = -4Z_{12}Z_{21},$$
implies that the zeros of $Z_{11} - Z_{22} - \sqrt{G}$ are included in the set of branch points $\mathcal{P}$.

We next define the $2 \times 2$-matrix valued function $\mathcal{B}(k)$ by
\begin{align}\label{Bdef}
\mathcal{B}(k) = \frac{\log z(k)}{\tau} S^b(k) \sigma_3 S^b(k)^{-1}, \qquad k \in \C \setminus \mathcal{C}.
\end{align}
By adding, if necessary, branch cuts to $\mathcal{C}$ to ensure that $\log z(k)$ is single valued on $\C \setminus \mathcal{C}$, 
we find $e^{\tau\mathcal{B}(k)} = Z(k)$ and
$$e^{-t\mathcal{B}(k)} = S^b(k) e^{it\tilde{\Omega}(k) \sigma_3} S^b(k)^{-1}, \qquad k \in \C \setminus \mathcal{C},$$
where 
$$\tilde{\Omega}(k) = -\frac{\log z(k)}{i\tau}.$$ 

Floquet theory implies that the matrix valued function $P(t,k)$ defined by
$$P(t,k) = \psi(t,k)  e^{-t\mathcal{B}(k)}, \qquad k \in \C\setminus \mathcal{C},$$
is $t$-periodic with period $\tau$. 
Defining $\mathcal{E}(t,k)$ by
\begin{align}\label{calEdef}
\mathcal{E}(t,k) = P(t,k) S^b(k), \qquad k \in \C \setminus \mathcal{C},
\end{align}
we deduce that the function $\psi^b(t,k)$ defined by
\begin{align}\label{psiB}
\psi^b(t,k) = \psi(t,k) S^b(k)
= \mathcal{E}(t,k) e^{-i\tilde{\Omega}(k) t \sigma_3}, \qquad k \in \C \setminus \mathcal{C},
\end{align}
is a solution of (\ref{tpartB}), where $\mathcal{E}(t,k)$ is time-periodic with period $\tau$.
We define the solution $\phi_1(x,t,k)$ of (\ref{lax}) by
$$\phi_1(x,t,k) = \mu_1(x,t,k) e^{-i(kx + \tilde{\Omega}(k) t)\sigma_3},$$
where $\mu_1(x,t,k)$ is the unique solution of the linear Volterra integral equation
\begin{align}\nonumber
\mu_1(&x,t,k) = e^{-ikx\hat{\sigma}_3} \bigg\{\mathcal{E}(t,k) 
-  \mathcal{E}(t,k) \int_t^\infty e^{i\tilde{\Omega}(k) (t' - t)\hat{\sigma}_3} \big[\mathcal{E}^{-1}(t', k) 
	\\\label{mu1def}
&\times (V - V^b)(0,t',k) \mu_1(0,t',k) \big] dt' 
+ \int_0^x e^{ikx' \hat{\sigma}_3}[U(x',t) \mu_1(x',t,k)] dx'\bigg\}.
\end{align}

Let \begin{align}\nonumber
D_1 = \{\im k > 0\} \cap \{\im \tilde{\Omega}(k) > 0\},  \qquad
D_2 = \{\im k > 0\} \cap \{\im \tilde{\Omega}(k) < 0\}, 
	\\ \label{Djdef}
D_3 = \{\im k < 0\} \cap \{\im \tilde{\Omega}(k) > 0\},  \qquad
D_4 = \{\im k < 0\} \cap \{\im \tilde{\Omega}(k) < 0\},
\end{align}
and let $D_+ = D_1 \cup D_3$ and $D_- = D_2 \cup D_4$.
The eigenfunctions $\{\mu_j(x,t,k)\}_1^3$ possess the following analyticity and boundedness properties:
\begin{itemize}
  \item The first (resp. second) column of $\mu_1(0,t,k)$ is defined and analytic for $D_- \setminus \mathcal{C}$ (resp. $D_+ \setminus \mathcal{C}$).
   
  \item $\mu_1(0,t,k)$ approaches $\mathcal{E}(t,k)$ as $t \to \infty$. More precisely, if $K_\pm$ are compact subsets of $(\overline{D_\pm \setminus \mathcal{C}})\setminus \mathcal{P}$, then
\begin{align}\label{mu1calE}
 |\mu_1(0,t,k) - \mathcal{E}(t,k)| \leq C(1+t)^{-5/2}, \qquad k \in (K_-, K_+), \quad t \geq 0,
\end{align}
 where the notation $k \in (K_-, K_+)$ indicates that the first and second columns are valid for $k \in K_-$ and $k \in K_+$, respectively. 
 
  \item $\mu_2(x,t,k)$ is defined and analytic for all $k \in \C$.

  \item The first (resp. second) column of $\mu_3(x,t,k)$ is defined and analytic for $\im k < 0$ (resp. $\im k > 0$) with a continuous extension to $\im k \leq 0$ (resp. $\im k \geq 0$).
\end{itemize}

\subsection{The global relation}
The expression
$$C(t, k) = e^{i\tilde{\Omega}(k)t\sigma_3}\mu_1^{-1}(0,t,k) \mu_3(0,t,k) e^{-2ik^2t\sigma_3}$$
is independent of $t$. Letting $t' \to \infty$ in the $(12)$ entry of the relation
\begin{align}\label{preGR}
C(t, k) = C(t',k),
\end{align}
we find
\begin{align}\nonumber
&(\mu_1(0,t,k))_{22} b(t,k) - (\mu_1(0,t,k))_{12} a(t,k) = 0, 
	\\ \label{gGR}
& \hspace{5cm} k \in D_1 \setminus \mathcal{C}, \quad \im(\tilde{\Omega}(k) + 2k^2) > 0, \quad t \geq 0,
\end{align}
where $a(t,k)$ and $b(t,k)$ are defined by
\begin{align}\label{atildebtildedef}
a(t,k) = (\mu_3(0,t,k))_{22}, \qquad b(t,k) = (\mu_3(0,t,k))_{12}.
\end{align}

\section{Proof of Theorem \ref{pertth}}\nequation\label{proofsec}
Before we present the detailed proof, we provide an argument motivating the relevant constructions. 

\subsection{Motivation}
Assuming that the branch cuts have been chosen in such a way that $D_1 \setminus \mathcal{C}$ is connected, equation (\ref{gGR}) is valid in all of $\overline{D_1 \setminus \mathcal{C}}$ by analytic continuation. Since $a(t,k)$ and $b(t,k)$ are continuous in $\im k \geq 0$, it follows that, away from the possible zeros of $a(t,k)$, the function
\begin{align}\label{pGR}
Q(t,k) = \frac{(\mu_1(0,t,k))_{12}}{(\mu_1(0,t,k))_{22}} = \frac{b(t,k)}{a(t,k)}, \qquad k \in D_1 \setminus \mathcal{C},
\end{align}
has a continuous extension to $\bar{D}_1$. 

On the other hand, equation (\ref{mu1calE}) shows that $Q(t,k)$ asymptotes towards the $t$-periodic function $Q^b(t,k)$ defined by
$$Q^b(t,k) = \frac{(\mathcal{E}(t,k))_{12}}{(\mathcal{E}(t,k))_{22}}$$
as $t \to \infty$, that is,
\begin{align}\label{ppB}
Q(t,k) \sim Q^b(t,k), \qquad t \to \infty, \quad k \in \bar{D}_1 \setminus \mathcal{C}.
\end{align}
Since $\psi^b(t,k) = \mathcal{E}(t,k) e^{-i\tilde{\Omega}(k) t \sigma_3}$ satisfies the background $t$-part, we find that $Q^b$ satisfies the Ricatti equation
\begin{align}\label{ratioeq}
& Q_t^b + \lambda(2k\bar{g}_0^b - i\bar{g}_1^b)(Q^b)^2 + (2i\lambda|g_0^b|^2 + 4ik^2)Q^b - (2kg_0^b + ig_1^b) = 0.
\end{align}
Denoting the Fourier coefficients of $g_0^b(t)$, $g_1^b(t)$, and $Q^b(t,k)$ by $a_n$, $c_n$, and $d_n(k)$ respectively,
\begin{align}
& g_0^b(t) = \sum_{n=-\infty}^\infty a_n e^{in\omega t}, \qquad
g_1^b(t) = \sum_{n=-\infty}^\infty c_n e^{in\omega t},
	\\
& Q^b(t,k) = \sum_{n=-\infty}^\infty d_n(k) e^{in\omega t}, \qquad t \to \infty, \quad k \in \bar{D}_1 \setminus \mathcal{C},
\end{align}
and substituting these Fourier series representations into the Ricatti equation (\ref{ratioeq}), we find
\begin{align*}
& \sum_{n=-\infty}^\infty \biggl\{
in\omega d_n(k) 
+ 2\lambda k \sum_{l,m=-\infty}^\infty \bar{a}_l d_m(k) d_{n+l-m}(k) 
 - \lambda i \sum_{l, m=-\infty}^\infty  \bar{c}_{l} d_m(k) d_{n+l-m}(k)
 	\\
& + 2i\lambda \sum_{l,m = -\infty}^\infty \bar{a}_l a_m d_{n+l-m}(k)
+ 4ik^2 d_n(k) - 2k a_n  - i c_n
\biggr\}e^{in\omega t} = 0, \qquad k \in \bar{D}_1 \setminus \mathcal{C}.
\end{align*}
This yields the following infinite hierarchy of algebraic equations:
\begin{align}\nonumber
d_n(k) = &\;\frac{1}{4ik^2 + in\omega}\biggl\{
-2\lambda k \sum_{l,m=-\infty}^\infty \bar{a}_l d_m(k) d_{n+l-m}(k) 
 + \lambda i \sum_{l, m=-\infty}^\infty  \bar{c}_{l} d_m(k) d_{n+l-m}(k)
	\\ \label{algebraicsystem}
&- 2i\lambda \sum_{l,m = -\infty}^\infty \bar{a}_l a_m d_{n+l-m}(k) + 2k a_n + i c_n\biggr\}, 
\quad n \in \Z, \quad k \in \bar{D}_1 \setminus \mathcal{C}.
\end{align}

At this stage we would like to: (a) Use the fact that $Q(t,k)$ is continuous in $\bar{D}_1$ together with equation (\ref{ppB}) and the time-periodicity of $Q^b(t,k)$ to conclude that $Q^b(t,k)$ (and hence also $d_n(k)$) is nonsingular in $\bar{D}_1$. (b) Use the nonsingularity of the $d_n(k)$'s to solve the algebraic system (\ref{algebraicsystem}) for the $c_n$'s and the $d_n(k)$'s. 

The nonlinear nature of (\ref{algebraicsystem}) together with the unknown structure of $\tilde{\Omega}(k)$, makes this a challenging problem in general. However, it is a remarkable fact that in the perturbative limit, the system (\ref{algebraicsystem}) can be solved uniquely to all orders. The key observations that make this perturbative solution possible can be summarized as follows:
\begin{enumerate}[1.]
  \item In the perturbative limit, the domains $\{D_j\}_1^4$ simplify to the four quadrants $\{D_j^0\}_1^4$.

\item In the perturbative limit, $a(t,k)$ has no zeros since $a(t,k) = 1 + O(\epsilon)$.

  \item In the perturbative limit, 
\begin{align*}
& G(k) = -4\sin^2(2k^2 \tau) + \epsilon^2 G_2(k) + \epsilon^3 G_3(k)\cdots, \qquad \epsilon \to 0,
	\\
& \mathcal{B}(k) = -2ik^2\sigma_3 + \epsilon \mathcal{B}_1(k) + \epsilon^2 \mathcal{B}_2(k) + \cdots, \qquad \epsilon \to 0,
	\\
& \mathcal{E}(t,k) = I + \epsilon \mathcal{E}_1(t,k) + \epsilon^2 \mathcal{E}_2(t,k) + \cdots, \qquad \epsilon \to 0,	
\end{align*}
where the coefficients $\{G_N(k)\}_2^\infty$, $\{\mathcal{B}_N(k)\}_1^\infty$, and $\{\mathcal{E}_N(t,k)\}_1^\infty$ are analytic in the whole complex $k$-plane except for possible poles at the zeros of $\sin(2k^2\tau)$. In particular, no branch cuts survive in the limit $\epsilon \to 0$.

\item The zeros of $\sin(2k^2\tau)$ occur at
$$k = \pm \frac{\sqrt{-n\omega}}{2}, \qquad n \in \Z,$$
which are also the points at which the denominator in (\ref{algebraicsystem}) vanishes. 
By enforcing that the poles of $d_n(k)$ at these points must be removable singularities, we can solve the system (\ref{algebraicsystem}) to all orders.
\end{enumerate}

\subsection{Proof of theorem \ref{pertth}}
Since we do not assume that the series in (\ref{uformal}) converges, the power series in $\epsilon$ utilized in this proof should be considered merely as a convenient way of organizing the coefficients of these series. We emphasize that the coefficients of these series are rigorously defined. 

{\bf Step 1.}
We define coefficients $\{\psi_N(t,k)\}_0^\infty$ of a perturbative solution $\psi = \sum_{N=0}^\infty \epsilon^N \psi_N$ of (\ref{tpartB}) satisfying $\psi(0,k) = I$ by substituting the expansions
\begin{align}\label{qggpert}
g_0^b(t) = \sum_{N=1}^\infty \epsilon^N g_{0N}^b(t),\qquad
g_1^b(t) = \sum_{N=1}^\infty \epsilon^N g_{1N}^b(t),
\end{align}
into (\ref{tpartB}) and solving the resulting equation order by order. This yields
\begin{align*}
 & \psi_0(t,k) = e^{-2ik^2t\sigma_3}, 
  	\\
 & \psi_1(t,k) = \begin{pmatrix} 0 & \int_0^t e^{4ik^2(t'-t)} (2 k g_{01}^b + i g_{11}^b) dt'  \\
 \lambda \int_0^t e^{-4ik^2(t'-t)} (2 k \bar{g}_{01}^b - i \bar{g}_{11}^b) dt'  & 0  \end{pmatrix}
 e^{-2ik^2t\sigma_3},
\end{align*}
etc. Defining entire functions $\{Z_N(k)\}_0^\infty$ by $Z_N(k) = \psi_N(\tau,k)$, it follows that $Z = \sum_{N=0}^\infty \epsilon^N Z_N$ satisfies
\begin{align*}
Z(k) & = e^{-2ik^2\tau\sigma_3}  
	\\
& + \epsilon  e^{-2ik^2\tau\sigma_3} \begin{pmatrix} 0 & \int_0^\tau e^{4ik^2t} (2 k g_{01}^b + i g_{11}^b) dt  \\
 \lambda  \int_0^\tau e^{-4ik^2t} (2 k \bar{g}_{01}^b - i \bar{g}_{11}^b) dt & 0 \end{pmatrix}
 + O(\epsilon^2).
 \end{align*}
We introduce coefficients $\{G_N(k)\}_1^N$ of $G(k)$ and $\{z_N(k)\}_1^N$ of $z(k)$ by substituting the above expression for $Z(k)$ into (\ref{lambdaGdef}):
\begin{align*}
& G(k) = -4\sin^2(2k^2 \tau) + O(\epsilon^2),
	\\
& \sqrt{G(k)} = 2i\sin(2k^2 \tau) \bigg[1 + O\bigg(\frac{\epsilon^2}{\sin^2(2k^2\tau)}\bigg) \bigg],
	\\
& z(k) = e^{-2ik^2\tau} + O\bigg(\frac{\epsilon^2}{\sin(2k^2\tau)}\bigg).
\end{align*}
We next define coefficients $\{S^b_N(k)\}_1^\infty$ of $S^b = I + \epsilon S^b_1 + \cdots$ by taking the $N$'th order terms of (\ref{mathcalSdef}). 
Each coefficient of $\sqrt{G}$, $z$, and $S^b$ can be written as an entire function divided by some power of $\sin(2k^2\tau)$. 
It follows that the same is true for the coefficients $\{\tilde{\Omega}_N(k)\}_0^\infty$, $\{\mathcal{B}_N(k)\}_0^\infty$, and $\{\mathcal{E}_N(t,k)\}_0^\infty$ defined by
\begin{align*}
& \tilde{\Omega}_N(k) = \frac{i}{\tau} \{\log z(k)\}_N, \qquad \mathcal{B}_N(k) = \bigg\{\frac{\log z(k)}{\tau} S^b(k) \sigma_3 S^b(k)^{-1}\bigg\}_N,
	\\
& \mathcal{E}_N(t,k) = \big\{\psi(t,k) e^{-t\mathcal{B}(k)} S^b(k) \big\}_N.
\end{align*}
In particular, these coefficients are analytic in the whole complex $k$-plane except for possible poles at the zeros of $\sin(2k^2\tau)$.
To leading order,
\begin{align*}
& \tilde{\Omega}(k) = 2k^2 + O(\epsilon^2), \qquad \mathcal{B}(k) = -2ik^2\sigma_3 + O(\epsilon), \qquad \mathcal{E}(t,k) = I + O(\epsilon).
\end{align*}

{\bf Step 2.}
We define functions $\mu_{jN}(x,t,k)$, $j = 1,2,3$, $N \geq 1$, by substituting the expansions
\begin{align}\label{mujpert}
\mu_j(x,t,k) = I + \epsilon \mu_{j1}(x,t,k) + \epsilon^2 \mu_{j2}(x,t,k) + \cdots, \qquad j = 1,2,3,
\end{align}
into the integral equations (\ref{mu23def}) and (\ref{mu1def}) and equating the coefficients of $\epsilon^N$ for each $N \geq 1$:
\begin{subequations}
\begin{align} \nonumber
 \mu_{1N}(x,t,k) = &\; e^{-ikx\hat{\sigma}_3} \bigg\{\mathcal{E}(t,k) 
- \mathcal{E}(t,k) \int_t^\infty e^{i\tilde{\Omega}(k) (t' - t)\hat{\sigma}_3} \big[\mathcal{E}^{-1}(t', k) 
	\\ \label{mu1Ndef}
&\times ((V - V^b) \mu_1)(0,t',k) \big] dt' 
+ \int_0^x e^{ikx' \hat{\sigma}_3}(U\mu_1)(x',t,k) dx'\bigg\}_N,
	\\\nonumber
 \mu_{jN}(x,t,k) = &\; \int_{(x_j, t_j)}^{(x,t)} e^{i[k(x'-x) + 2k^2(t'-t)] \hat{\sigma}_3} \{(U\mu_j)(x',t',k) dx' + (V \mu_j)(x',t',k) dt'\}_N, 
 	\\ \label{mujNdef}
&\hspace{8cm}  j = 2,3.
\end{align}
\end{subequations}
These equations define the coefficients $\mu_{jN}(x,t,k)$ recursively for all $N \geq 1$. For $N =1$, we have
\begin{align*}
 \mu_{11}(x,t,k) = &\; e^{-ikx\hat{\sigma}_3} \bigg[\mathcal{E}_1(t,k) 
- \int_t^\infty e^{2ik^2 (t' - t)\hat{\sigma}_3} (V_1 - V_1^b)(0,t',k)  dt' 
	\\
& + \int_0^x e^{ikx' \hat{\sigma}_3} U_1(x',t) dx'\bigg],
	\\
\mu_{j1}(x,t,k) = &\; \int_{(x_j, t_j)}^{(x,t)} e^{i[k(x'-x) + 2k^2(t'-t)] \hat{\sigma}_3} [U_1(x',t') dx' + V_1(x',t',k) dt'], \qquad j = 2,3,
\end{align*}
where
$$U_1 = \begin{pmatrix} 0 & u_1 \\
\lambda \bar{u}_1 & 0 \end{pmatrix}, \qquad V_1 = \begin{pmatrix} 0 & 2ku_1 + iu_{1x} \\
2\lambda k \bar{u}_1 - i\lambda \bar{u}_{1x} & 0 \end{pmatrix}.$$
Using the expansion 
\begin{align}\label{eitildeOmega}
e^{i\tilde{\Omega}(k) t} = e^{2 ik^2 t} \big(1 + \epsilon^2 it\tilde{\Omega}_2(k) + O(\epsilon^3)\big)
\end{align}
and the assumption that $V_1 - V_1^b = O(t^{-7/2})$, we deduce that the $\mu_{jN}$'s have the following properties for any $N \geq 1$.
\begin{itemize}
  \item The first (resp. second) column of $\mu_{1N}(x,t,k)$ is defined and analytic for $k \in D_2^0$ (resp. $k \in D_3^0$). For $x = 0$, the first (resp. second) column of $\mu_{1N}(0,t,k)$ is defined and analytic for $k \in D_-^0$ (resp. $k \in D_+^0$).
   
  \item $\mu_{2N}(x,t,k)$ is defined and analytic for all $k \in \C$.

  \item The first (resp. second) column of $\mu_{3N}(x,t,k)$ is defined and analytic for $\im k < 0$ (resp. $\im k > 0$) with a continuous extension to $\im k \leq 0$ (resp. $\im k \geq 0$).
 
  \item As $k \to \infty$, $\mu_{1N}(0, t, k) = O(k^{-1})$ where the first column is valid for $\arg k \in (\frac{\pi}{2} + \delta, \pi - \delta) \cup (\frac{3\pi}{2} + \delta, 2\pi - \delta)$ and the second column is valid for $\arg k \in (\delta, \frac{\pi}{2}  - \delta) \cup (\pi + \delta, \frac{3\pi}{2}  - \delta)$, where $\delta > 0$ is arbitrarily small. Moreover, as $k \to \infty$,
   \begin{subequations}
   \begin{align}
 &  \mu_{2N}(0, t, k) = O(k^{-1}), \qquad k \in (\bar{D}_+^0, \bar{D}_-^0),
   	\\
&   \mu_{3N}(x, t, k) = O(k^{-1}),  \qquad k \in (\bar{\C}_-, \bar{\C}_+).
   \end{align}
   \end{subequations}
\end{itemize}

{\bf Step 3.}
Equation (\ref{pGR}) suggests defining functions $\{Q_N(t,k)\}_1^\infty$ by
\begin{align*}
Q_N(t,k) & = \bigg\{\frac{(\mu_1(0,t,k))_{12}}{(\mu_1(0,t,k))_{22}}\bigg\}_N 
	\\
& = (\mu_{1N}(0,t,k))_{12}
- (\mu_{1,N-1}(0,t,k))_{12}(\mu_{11}(0,t,k))_{22} + \cdots, \qquad k \in D_+^0.
\end{align*}
The terms of order $O(\epsilon^N)$ of equation (\ref{preGR}) imply
\begin{align}\label{preGR2}
 C_N(t, k) = C_N(t',k),
\end{align}
where
$$C_N(t,k) = \bigg\{e^{i\tilde{\Omega}(k)t\sigma_3}\mu_1^{-1}(0,t,k) \mu_3(0,t,k) e^{-2ik^2t\sigma_3}\bigg\}_N.$$
Letting $t' \to \infty$ in the $(12)$ entry of (\ref{preGR2}), using (\ref{eitildeOmega}) and the fact that the $\mu_{1N}$'s and $\mu_{3N}$'s are defined and continuous for $k \in D_1^0$, we find 
\begin{align}\label{gGR2}
\Big\{e^{i\tilde{\Omega}(k)t}[(\mu_1(0,t,k))_{22} b(t,k) - (\mu_1(0,t,k))_{12} a(t,k) ]\Big\}_N = 0, \quad k \in D_1^0, \;\; t \geq 0,
\end{align}
where $a$ and $b$ are defined in (\ref{atildebtildedef}).
Since equation (\ref{gGR2}) is valid for all $N \geq 1$, it follows that
\begin{align*}
\Big\{(\mu_1(0,t,k))_{22} b(t,k) - (\mu_1(0,t,k))_{12} a(t,k)\Big\}_N = 0, \qquad k \in D_1^0, \quad t \geq 0.
\end{align*}
Hence
\begin{align}\label{QNeq}
Q_N(t,k) = \bigg\{\frac{b(t,k)}{a(t,k)}\bigg\}_N
= b_N(t,k) - b_{N-1}(t,k) a_1(t,k) + \cdots, \qquad k \in D_1^0,
\end{align}
where
$$b(t,k) =  \sum_{N=1}^\infty \epsilon^N b_N(t,k), \qquad 
a(t,k) = 1 + \sum_{N=1}^\infty \epsilon^N a_N(t,k).$$
In particular, $Q_N(t,k)$ is a continuous function of $k \in \bar{D}_1^0$ for all $N \geq 1$ and $t \geq 0$.

{\bf Step 4.}
Define $\{Q_N^b(t,k)\}_1^\infty$ by
$$Q_N^b(t,k) = \bigg\{ \frac{(\mathcal{E}(t,k))_{12}}{(\mathcal{E}(t,k))_{22}}\bigg\}_N
= (\mathcal{E}_N(t,k))_{12} - (\mathcal{E}_{N-1}(t,k))_{12}(\mathcal{E}_1(t,k))_{22} + \cdots, \qquad k \in \C.$$
The analyticity properties of the $\mathcal{E}_N$'s imply that the functions $Q_N^b(t,k)$ are analytic for $k \in \C$ except for possible poles at the zeros of $\sin(2k^2\tau)$.
We claim that $Q_N^b(t,k)$ cannot have poles at the zeros of $\sin(2k^2\tau)$ in $\bar{D}_1^0$. 
To see this, note that the second column of (\ref{mujNdef}) yields
\begin{align*}
\begin{cases}
(\mu_{3N}(x,t,k))_{12} = - \int_x^\infty e^{2ik(x'-x)}  \Big\{ u(x',t) (\mu_3(x',t,k))_{22}  \Big\}_N dx',
	\\
(\mu_{3N}(x,t,k))_{22} = - \int_x^\infty \lambda \Big\{ \bar{u}(x',t) (\mu_3(x',t,k))_{12}\Big\}_N dx',
\end{cases}\qquad N \geq 1.
\end{align*}
Together with the assumption that the set $\{u_N(\cdot, t)| t \geq 0\}$ is bounded in $L^1([0,\infty))$ this implies  
$$\sup_{\substack{\im k \geq 0 \\ t \geq 0}} |a_N(t,k)| < \infty, \qquad 
\sup_{\substack{\im k \geq 0 \\ t \geq 0}} |b_N(t,k)| < \infty, \qquad N \geq 1.$$
In particular, by (\ref{QNeq}), there exist constants $C_N > 0$ such that
\begin{align}\label{pNbounded}
\sup_{\substack{\im k \geq 0 \\ t \geq 0}} Q_N(t, k) \leq C_N, \qquad N \geq 1.
\end{align}
Let $K \in \bar{D}_1^0$ be a zero of $\sin(2k^2\tau)$.
Fix $N \geq 1$ and $t_0 \geq 0$. 
Equation (\ref{mu1Ndef}) implies the pointwise convergence
\begin{align}\label{ppNB}
\lim_{t \to \infty} |Q_N(t,k) - Q_N^b(t,k)| = 0, \qquad k \in D_1^0, \quad N \geq 1.
\end{align}
For each $k \in D_1^0$, (\ref{ppNB}) shows that there exists a $T_k> 0$, such that $|Q_N(t,k) - Q_N^b(t,k)| < 1$ whenever $t > T_k$. For each $k \in D_1^0$, let $n_k$ be an integer such that $t_0 + n_k \tau > T_k$. The periodicity of $Q_N^b$ then implies
\begin{align*}
& |Q_N^b(t_0,k)| = |Q_N^b(t_0 + n_k \tau,k)| 
	\\
& \leq |Q_N^b(t_0 + n_k \tau, k) - Q_N(t_0 + n_k \tau, k)| + |Q_N(t_0 + n_k \tau, k)|
\leq 1 + C_N,
\end{align*}
showing that $Q_N^b(t_0,k)$ is a bounded function of $k \in D_1^0$. It follows that $Q_N^b$ cannot have any poles in $\bar{D}_1^0$ and is in fact continuous for $k \in \bar{D}_1^0$.

{\bf Step 5.}
Let $a_{N,n}$, $c_{N,n}$ denote the Fourier coefficients of $g_{0N}^b(t)$, $g_{1N}^b(t)$ as in (\ref{g0NBFourier}) and (\ref{g1NBFourier}). Moreover, let 
$$d_{N,n}(k) = \frac{1}{\tau} \int_0^\tau Q_N^b(t,k) e^{-in\omega t} dt,$$
denote the Fourier coefficients of the $t$-periodic function $Q_N^b(t,k)$:
\begin{align}\label{pNBFourier}
& Q_N^b(t,k) = \sum_{n=-\infty}^\infty d_{N,n}(k) e^{in\omega t}, \qquad t \to \infty, \quad N \geq 0, \quad k \in \bar{D}_1^0.
\end{align}
Since $Q_N^b$ is continuous for $k \in \bar{D}_1^0$ so is $d_{N, n}(k)$ for each $N \geq 1$ and $n \in \Z$.
The Ricatti equation (\ref{ratioeq}) is satisfied to all orders in $\epsilon$; the terms of order $O(\epsilon^N)$ imply that the functions $\{Q_N^b(t,k)\}_1^\infty$ satisfy the following hierarchy of equations:
\begin{align}\nonumber
& \Big\{Q_{t}^b + \lambda(2k\bar{g}_0^b - i\bar{g}_1^b)(Q^b)^2 + (2i\lambda|g_0^b|^2 + 4ik^2)Q^b - (2kg_0^b + ig_1^b) \Big\}_N = 0, 
	\\ \label{Ricattipert}
& \hspace{7cm} k \in \bar{D}_1^0, \quad N \geq 1.
\end{align}
Substituting the Fourier series (\ref{g0NBFourier}), (\ref{g1NBFourier}), and (\ref{pNBFourier}) into (\ref{Ricattipert}), we find
\begin{align*}
& \sum_{n=-\infty}^\infty \biggl\{
in\omega d_n(k) 
+ 2\lambda k \sum_{l,m=-\infty}^\infty \bar{a}_l d_m(k) d_{n+l-m}(k) 
 - \lambda i \sum_{l, m=-\infty}^\infty  \bar{c}_{l} d_m(k) d_{n+l-m}(k)
 	\\
& + 2i\lambda \sum_{l,m = -\infty}^\infty \bar{a}_l a_m d_{n+l-m}(k)
+ 4ik^2 d_n(k) - 2k a_n  - i c_n
\biggr\}_N e^{in\omega t} = 0, \qquad k \in \bar{D}_1^0.
\end{align*}
This yields the infinite hierarchy of algebraic equations (\ref{dNndef}). The requirement that $d_{N,n}(k)$ has a removable singularity at each zero $k_1(n)$ of $4k^2 + n\omega = 0$ in $\partial D_1^0$ then yields (\ref{cNndef}).
This completes the proof.
\proofend

\begin{remark}\upshape\label{L1remark}
In the above proof, we assumed that $\|u_N(\cdot, t)\|_{L^1([0,\infty))}$ remains bounded as $t \to \infty$ for each $N \geq 1$. This assumption can be weakened if $V$ approaches $V^b$ faster than $O(t^{-7/2})$.
For example, suppose (\ref{uNgN}) is replaced with
\begin{align}\label{uNx}
u_{N}(0, \cdot) - g_{0N}^b \in \mathcal{S}([0,\infty)), \quad
u_{Nx}(0, \cdot) - g_{1N}^b \in \mathcal{S}([0,\infty)), \qquad N \geq 1.
\end{align}
Then the conclusion of Theorem \ref{pertth} remains valid under the weaker assumption that, for each $N \geq 1$,
\begin{align}\label{growsslower}
\text{$\|u_N(\cdot, t)\|_{L^1([0,\infty))}$ grows slower than some power of $t$ as $t \to \infty$.}
\end{align}
Indeed, fix $N \geq 1$ and let $K \in \bar{D}_1^0$ be a zero of $\sin(2k^2\tau)$.
By (\ref{growsslower}), there exist a $P > 0$ such that (cf. equation (\ref{pNbounded}))
\begin{align}\label{pNgrowth}
  Q_N(t, k) \leq C t^{P}, \qquad \im k \geq 0, \quad t \geq 1.
\end{align}
On the other hand, it follows from (\ref{uNx}) and (\ref{mu1Ndef}) that there exists a $B \geq 0$ such that, for any $A > 0$,
\begin{align}\label{pNpNB}
|Q_N(t,k) - Q_N^b(t,k)| \leq \frac{C}{t^A |k - K|^B},
\end{align}
for all $k \in D_1^0$ close to $K$ and all $t \geq 1$. Let $t_0 \geq 1$. Using the periodicity of $Q_N^b$, the inequality (\ref{pNgrowth}), and the estimate (\ref{pNpNB}) with  $A = 4 B P$, we find
\begin{align*}
|Q_N^b(t_0,k)| & = |Q_N^b(t_0 + n \tau,k)| 
	\\
& \leq |Q_N^b(t_0 + n \tau, k) - Q_N(t_0 + n \tau, k)| + |Q_N(t_0 + n \tau, k)|
 	\\
& \leq \frac{C}{(t_0 + n \tau)^{4BP} |k - K|^B} + C (t_0 + n \tau)^{P},
\end{align*}
for each integer $n \geq 0$. Now assume $k(s)$, $s \in [0,1]$, is a continuous path in $D_1^0$ such that $k(0) = K$. 
For each sufficiently small $s$, there exists an $n = n(s)$ such that 
$$|k(s) - K|^{-\frac{1}{4P}} \leq t_0 + n(s) \tau \leq |k(s) - K|^{-\frac{1}{2P}}.$$
Then $(t_0 + n(s) \tau)^{-4BP} |k(s) - K|^{-B} \leq 1$ and $ (t_0 + n(s) \tau)^{P} \leq |k(s) -K|^{-1/2}$. Hence
$$|Q_N^b(t_0,k(s))| \leq C + C |k(s) - K|^{-1/2},$$
for all small enough Ê$s$. This shows that $Q_N^b(t_0, \cdot)$ cannot have a pole at $K$. Thus $Q_N^b(t_0,k)$ is a bounded function of $k \in D_1^0$ and the conclusion of Theorem \ref{pertth}  follows as in Step 5 above.
\end{remark}

\begin{remark}\upshape
The above proof shows that the perturbative coefficients $\{Q_N(t,k)\}_1^\infty$ and $\{Q_N^b(t,k)\}_1^\infty$ do not have poles in $\partial D_1^0$. This result only holds in the perturbative limit. Indeed, consider the example of stationary one-solitons. For these solitons (see \cite{tperiodicI})
$$Q(t,k) = \frac{(\mu_1(0,t,k))_{12}}{(\mu_1(0,t,k))_{22}}
= \frac{\sqrt{\omega} e^{it\omega}}{\sqrt{\omega} \sinh(\gamma) + 2ik \cosh (\gamma )}$$
and
$$Q^b(t,k) = \frac{(\mathcal{E}(t,k))_{12}}{(\mathcal{E}(t,k))_{22}} 
= \frac{\sqrt{\omega } e^{i t \omega }}{\sqrt{\omega } \sinh (\gamma )+2 i k \cosh (\gamma )}.$$
Hence $Q(t,k) = Q^b(t,k)$ has a simple pole at $k_0 = \frac{i}{2} \sqrt{\omega } \tanh (\gamma )$.
If $\gamma \geq 0$, $k_0$ lies in $\partial D_1^0$. The existence of this pole is related to the zero of $a(t,k)$ which generates the soliton. Indeed, 
$$a(t,k) = \frac{2 k-i \sqrt{\omega } \tanh(\gamma )}{2 k+i \sqrt{\omega }}$$
has a simple zero at $k_0$. Thus, for $\gamma \geq 0$, both sides of the generalized global relation (\ref{pGR}) have simple poles in $\partial D_1^0$, whereas in the perturbative limit no such pole may exist.
\end{remark}

\begin{remark}\upshape
The perturbative expansion of the exponential $e^{i\tilde{\Omega}(k) t}$ typically generates secular terms involving powers of $t$, see (\ref{eitildeOmega}). One advantage of considering the quotient $Q(t,k)$, which does not involve $e^{i\tilde{\Omega}(k) t}$, is that these secular terms are avoided. 
\end{remark}

\bigskip
\noindent
{\bf Ethics statement.} This work did not involve any collection of human data.

\medskip
\noindent
{\bf Data accessibility.} This work does not rely on any experimental data.

\medskip
\noindent
{\bf Competing interests.} We have no competing interests.

\medskip
\noindent
{\bf Authors' contributions.} JL and ASF conceived the ideas and proved the mathematical results presented in this paper. Both authors gave final approval for publication.

\medskip
\noindent
{\bf Acknowledgements.} The authors are grateful to the two referees for many helpful suggestions.

\medskip
\noindent
{\bf Funding statement.} This work was supported by the EPSRC, grant EP/H04261X/1.

\bibliographystyle{plain}
\bibliography{is}

\end{document}